\documentclass[12pt]{article}
\usepackage{mathrsfs}
\usepackage{amsfonts}
\usepackage{amsmath,amssymb,amsfonts,amsthm,fancyhdr,mathrsfs}
\usepackage{bbm}
\usepackage{verbatim}
\usepackage{color}
\usepackage{dsfont}
\usepackage[left=2.5cm,right=2.5cm,top=3cm,bottom=2.5cm]{geometry}

\usepackage{latexsym}
\usepackage{graphicx}
\usepackage{enumerate}
\usepackage[colorlinks=true,urlcolor=blue,
citecolor=red,linkcolor=blue,linktocpage,pdfpagelabels,
bookmarksnumbered,bookmarksopen]{hyperref}
\usepackage{cite}
\usepackage{indentfirst}
\usepackage{amsmath}
\usepackage{booktabs}
\usepackage{threeparttable}
\usepackage{float}
\allowdisplaybreaks[4]
\numberwithin{equation}{section} 

\numberwithin{equation}{section}
\newtheorem{thm}{Theorem}[section]

\newtheorem{lem}[thm]{Lemma}

\newtheorem{re}{Remark}[section]
\newenvironment{pf}{{\noindent \it \bf Proof:}}{{\hfill $\square$}\\}

\title{Optimal convergence rate to the nonrelativistic limit of  Chandrasekhar variational model for Neutron stars}

\author{
	{\bf\large  Yuanhui Chen, 
		\bf\large Qingxuan Wang}\thanks{Qingxuan Wang is the corresponding author.}
	\thanks{E-mail:chyhhui@163.com(Y. Chen), wangqx@zjnu.edu.cn(Q. Wang).}\\
	\small School of Mathematical Sciences,
	{Zhejiang Normal University, Jinhua, 321004, China}\\
}
\date{}

\begin{document}
	\maketitle
	\begin{abstract}
In this paper, we consider the nonrelativistic limit of Chandrasekhar variational model for neutron stars.  We prove that  the minimizer $\rho_{c}$ of Chandrasekhar energy $E_c(N)$ converges strongly to   the minimizer $\rho_{\infty}$ of limit energy $E_{\infty}(N)$
in  $L^1\cap L^{\frac{5}{3}}(\mathbb{R}^3)$ as the speed of light $c\rightarrow\infty$, this is a limit between two free boundary problems. 
We develop a novel approach to obtain  that the above nonrelativistic limit has the optimal convergence rate $\frac{1}{c^2}$. For the radius $R_c$ of the compact support of $\rho_c(x)$ and the radius $R_\infty$ of the compact support of $\rho_\infty(x)$, we also get the optimal convergence rate, $R_\infty-R_c=O(\frac{1}{c^2})$ as $c\rightarrow\infty$. Moreover, we  obtain the optimal uniform bounds of $R_c$ and  $L^\infty$-norm of $\rho_c$   with respect to $N$ as $c\rightarrow \infty$.

%Moreover, we show that  the $L^\infty$-convergence rate of $|\rho_c-\rho_\infty|$ is not bigger than  $O(\frac{1}{c^3})$ in the corner layer $B(R_\infty +\frac{K_1}{c^2})\setminus B(R_\infty -\frac{K_1}{c^2})\, (K_1>0)$, while $L^\infty$-convergence rate of $|\rho_c-\rho_\infty|$ is not bigger than  $O(\frac{1}{c^2})$  inside the  corner layer $B(R_\infty -\frac{K_1}{c^2})$ for $N$ large enough.

	\end{abstract}
\textbf{Keywords}: Free boundary problem; Nonrelativistic limit; Optimal convergence rate\\
	\textbf{Mathematics Subject Classification (2020): 35Q75; 49S05; 49J45; 46T20}  
\section{Introduction and main results}\label{sec1}

\subsection{Background and Motivations}
A neutron star is a system of identical relativistic fermions which interact through the self-gravitational force. It is  well known that a neutron star collapses when its mass is bigger than a critical number. The maximum mass of a stable star, called the Chandrasekhar limit, which was computed by Chandrasekhar \cite{AJ-maximum-mass-dwarf}. The rigorous derivation of the Chandrasekhar functional  from many-body quantum theory has been done by Lieb and Yau in \cite{CMP-lieb-Chandrasekhar-theory}.
The Chandrasekhar theory is the relativistic analog of the famous Thomas-Fermi theory of nonrelativistic electrons in atomic physics \cite{ARMAgf}. In the Chandrasekhar theory, the ground state energy of a neutron (Fermion) star is given by
\begin{align}\label{He}
	E^{Ch}(N)=\inf\left\lbrace \mathcal{E}^{Ch}(\rho):0\leq\rho\in L^1\cap L^{4/3}(\mathbb{R}^3) \text{ and }\int_{\mathbb{R}^3}\rho(x)\,dx=N\right\rbrace
\end{align}
with the energy functional
\begin{gather}\label{Hfunctional}
	\begin{aligned} 
		\mathcal{E}^{Ch}(\rho)=\int_{\mathbb{R}^3}j_m(\rho(x))\,dx
		-\frac{\kappa}{2}\int_{\mathbb{R}^3}\int_{\mathbb{R}^3}\frac{\rho(x)\rho(y)}{|x-y|}\,dxdy,
	\end{aligned}
\end{gather}
where $N>0$ is the particle number, $m>0$ is the neutron mass and $\kappa/|x|$ represents the Newtonian gravitational potential in appropriate physical units. Here $\rho(x)$ is the density of the relativistic system and the  functional $j_m(\rho(x))$ is the semiclassical approximation for the relativistic kinetic energy  at density $\rho(x)$, namely
\begin{gather*}
	\begin{aligned} 
		j_m(\rho(x))&=\frac{q}{8\pi^3}\int_{|p|<\eta}\left( \sqrt{|p|^2+m^2}-m\right) \,dp\\
		&=\frac{q}{16\pi^2}\left[ \eta(2\eta^2+m^2)\sqrt{\eta^2+m^2}-m^4\ln\left( \frac{\eta+\sqrt{\eta^2+m^2}}{m}\right) \right]-m\rho,
	\end{aligned}
\end{gather*}
where $p\in \mathbb{R}^3$, $\sqrt{|p|^2+m^2}$ denotes a pseudo-relativistic kinetic operator, $\eta=(6\pi^2\rho/q)^{\frac{1}{3}}$ and for generality we assume the permutations of space-spin variables $q\ge1$, $q=2$ in nature, $q=1$ would correspond to spin-polarized matter \cite{CMP-lieb-Chandrasekhar-theory}. 
Recently, Nguyen\cite{JMPNg} studies blow-up profile  for the Chandrasekhar problem \eqref{He} as $\tau=\kappa N^{2/3}\uparrow \tau_c$ where $\tau_c$ is a critical constant, and then  as a continuation, they also obtained the similar blow-up profiles of neutron stars in the Hartree-
 Fock-Bogoliubov (HFB) theory\cite{CVNg}. The more relevant information and results of this model, we can refer to \cite{JFALopes,Lez-M,JMPHF}.

 \textbf{Our first motivation} of this paper comes from the work of Lenzmann\cite{Lenz-uniqueness}, Guo and Zeng\cite{Guo-Zeng}, who studied the uniqueness (referred to as the Lieb-Yau conjecture in \cite{CMP-lieb-Chandrasekhar-theory}) of the ground state for pseudo-relativistic Hartree energy using the non-relativistic limit of pseudo-relativistic Hartree energy, where the  functional is given by 
\begin{align*}
	\mathcal{E}^H_{c}(\psi)=\int_{\mathbb{R}^3 } \bar{\psi } (\sqrt{-c^2\Delta +m^2c^4}-mc^2) \psi\,dx-\frac{1}{2} \int_{\mathbb{R}^3 } (|x|^{-1}*|\psi |^2) |\psi |^2\,dx
\end{align*}
where $c$ denotes the speed of light. There are other physical settings considering 
 nonrelativistic limit, one can see \cite{Cvguo,2016-nonrelativistic,NL-of-Dirac-Fock-eq,NL-of-KGM,NL-of-nonlinear-Dirac-equation,25-NL-of-NKG,NLD} and reference therein.  In the general physics, the pseudo-relativistic kinetic operator is  given by $\sqrt{c^2|p|^2+m^2c^4}$ depending on the the speed of light $c$.  In this paper, we reinstall the speed of light $c>0$ into $\mathcal{E}^{Ch}(\rho)$ defined in (\ref{Hfunctional}), which yields the $c$-depending Chandrasekhar functional 	$\mathcal{E}_c(\rho)$. More precisely, we consider the following Chandrasekhar variational minimization problem
 \begin{align}\label{E-N}
 	E_c(N)=\inf\{\mathcal{E}_c(\rho):0\le\rho\in L^1\cap L^{\frac{4}{3}}(\mathbb{R}^3) \text{ and } \int_{\mathbb{R}^3}\rho(x)\,dx=N\},
 \end{align} 
 where the energy functional $\mathcal{E}_c(\rho)$ is given by 
 \begin{align}\label{eq1.5}
 	\mathcal{E}_c(\rho(x))=\int_{\mathbb{R}^3}j_{mc}(\rho(x))\,dx-\kappa D(\rho,\rho),
 \end{align}
 $j_{mc}(\rho(x))$ is given by 
 \begin{gather}\label{jmc}
 	\begin{aligned}
 		j_{mc}(\rho(x))&=\frac{q}{8\pi^3}\int_{|p|<(6\pi^2\rho/q)^{\frac{1}{3}}}\left( \sqrt{c^2|p|^2+m^2c^4}-mc^2\right) \,dp,
 	\end{aligned}
 \end{gather}
 and $D(\rho,\rho)$ is given by 
 \begin{gather}
 	\begin{aligned} 
 		D(\rho,\rho)=\frac{1}{2}\int_{\mathbb{R}^3}\int_{\mathbb{R}^3}\frac{\rho(x)\rho(y)}{|x-y|}\,dxdy.
 	\end{aligned}
 \end{gather}
 We mainly focus on the limiting behavior of minimizers of \eqref{E-N} as  $c\to\infty$,  which is called the \textbf{nonrelativistic limit}, see Theorem \ref{th1.2} below. 
 
 \textbf{Our second motivation} comes from the work  about the Ginzburg-Landau energy (see \cite{Brezis1993})
 \begin{gather*}
 	\begin{aligned}
 		E^{GL}_\varepsilon(u)=\frac{1}{2}\int_{\Omega}|\nabla u|^2\,dx
 		+\frac{1}{4\varepsilon^2}\int_{\Omega}\left(|u|^2-1\right)^2\,dx
 	\end{aligned}
 \end{gather*}
 minimized in $\mathscr{H}_g=\{u\in H^1(\Omega): \int_{\Omega}|u|^2\,dx=1, u=g(x)\ \ \text{in}\ \ \partial \Omega\}$
 and  the work  about the Gross-Pitaevskii energy (see \cite{ARMA2015} and references therein)
 \begin{gather*}
 	\begin{aligned}
 		E^{GP}_\varepsilon(u)=\frac{1}{2}\int_{\mathbb{R}^2}|\nabla u|^2\,dx+\frac{1}{2\varepsilon^2}\int_{\mathbb{R}^2}V(x)|u|^2\,dx
 		+\frac{1}{4\varepsilon^2}\int_{\mathbb{R}^2}|u|^4\,dx
 	\end{aligned}
 \end{gather*}
 minimized in $\mathscr{H}_N=\{u\in H^1(\mathbb{R}^2): \int_{\mathbb{R}^2} V(x)|u |^2<\infty, \int_{\mathbb{R}^2}|u|^2\,dx=N\}$.  They studied the limit behavior as $\varepsilon\rightarrow 0$, and can obtain $\|u_\varepsilon-u_0\|_{L^{\infty}(\Omega)}\leq{O}(\varepsilon^2)$ as $\varepsilon\to 0$ (may be called as \textit{Thomas-Fermi limit}), where $u_\varepsilon$ is a minimizer for $E^{GL}_\varepsilon(\cdot)$ (or $E^{GP}_\varepsilon(\cdot)$) on $\mathcal{H}_g$ and $u_0$ is identified. It is natural to ask whether the nonrelativistic limit mentioned above may obtain similar convergence rates. We should point out that these two limits are very different. The Thomas-Fermi limit is a limit from an elliptic variational problem to a free boundary problem (or Thomas-Fermi type problem), whereas the nonrelativistic limit in this paper is a limit between two free boundary problems. We develop a novel approach to obtain the convergence rates, see Theorem \ref{th1.3} below.
 
\subsection{Main results}
 Our first result is to show the minimizer $\rho_c$ of problem \eqref{E-N} converge strongly in $L^1\cap L^{\frac{5}{3}}(\mathbb{R}^3)$ to a minimizer $\rho_\infty$ of the following  minimizing problem
	\begin{align}\label{eq-limite-minimizing-problem}
	E_\infty(N)=\inf\{\mathcal{E}_\infty(\rho):\rho\in  L^1\cap L^{\frac{5}{3}}(\mathbb{R}^3) \text{ and }\int_{\mathbb{R}^3}\rho(x)\,dx=N\},
\end{align}
where $\mathcal{E}_\infty(\rho)$ is given by
\begin{gather}\label{eq-limiting-eneygy-functional}
	\begin{aligned}
		\mathcal{E}_\infty(\rho(x))=\int_{\mathbb{R}^3}j_{\infty}(\rho(x))\,dx-\kappa D(\rho,\rho),
	\end{aligned}
\end{gather}
here $j_{\infty}(\rho(x))$ is given by
\begin{align}\label{j-infty-rho}
	j_{\infty}(\rho(x))=\frac{q}{8\pi^3}\int_{|p|<\eta}\frac{|p|^2}{2m}\,dp=\frac{3}{10m}\big(\frac{6\pi^2}{q}\big)^{\frac{2}{3}}\rho(x)^{\frac{5}{3}},
\end{align}
where $\eta=(6\pi^2\rho/q)^{\frac{1}{3}}$. Note that  the operator $|p|^2$ is the classical kinetic operator, that is why we call the above limit as \textbf{nonrelativistic limit}.

  Before studying  the   nonrelativistic limit, we need to know the existence of the minimizers of $E_c(N)$ and $E_\infty(N)$. 
In the Chandrasekhar theory, the stellar collapse of big neutron stars boils down to the fact that $E^{Ch}(N)=-\infty$ if $N>N_*$. Let $ N_*$ be the  constant related to the Chandrasekhar limit $\kappa N_*^{2/3}$ in 
\cite{CMP-lieb-Chandrasekhar-theory},  or  the optimal constant of the following inequality
\begin{gather}\label{GNS2}
	\begin{aligned} 
		\kappa N_*^{2/3}D(\rho,\rho)
		\leq K_{cl}\|\rho\|_{L^{\frac{4}{3}}}^{\frac{4}{3}}\|\rho\|_{L^{1}}^{\frac{2}{3}},\quad \forall \ 0\leq\rho\in L^1\cap L^{4/3}(\mathbb{R}^3).
	\end{aligned}
\end{gather}
where $K_{cl}=\frac{3}{4}\left({6\pi^2}/{q}\right)^{1/3}$. We have the following existence theorem.
\begin{thm}\label{th1.1}
	Suppose $m>0$. Then we have \\
 (i) $E_c(N)$ has at least one  symmetric decreasing minimizer with compact support if and only if $0<N^{\frac{2}{3}}<cN_*^{\frac{2}{3}}$.\\
 (ii)  $E_\infty(N)$ has at least one  symmetric decreasing minimizer with compact support  for all $N>0$.
\end{thm}
\begin{re}
The existence of minimizer of   $E_\infty(N)$ is given in Appendix. For $E_c(N)$,  define 
\begin{align}\label{eq-minimizer-relation}
	\rho(x)=c^{-3}\tilde{\rho}(c^{-\frac{1}{2}}x),
\end{align}
an elementary calculation shows that
\begin{align}\label{eq1.17}
	\mathcal{E}^{Ch}(\rho)=c^{-\frac{7}{2}}\mathcal{E}_c(\tilde{\rho})\ \ \text{and} \int_{\mathbb{R}^3}\tilde{\rho}(x)\,dx=c^{\frac{3}{2}}\int_{\mathbb{R}^3}\rho(x)\,dx.
\end{align}
Thanks to the work of Lieb and Yau in 
\cite{CMP-lieb-Chandrasekhar-theory}, we know that for each $0<N<N_*$, there exists a symmetric decreasing minimizer $\rho^{Ch}$ with compact support for $E^{Ch}(N)$. Therefore, $\tilde{\rho}$ minimizes $\mathcal{E}_c({\rho})$ subject to $\int_{\mathbb{R}^3}\rho(x)\,dx=N$ if and only if $	\rho=c^{-3}\tilde{\rho}(c^{-\frac{1}{2}}x)$ minimizes $\mathcal{E}^{Ch}(\rho)$ subject to $\int_{\mathbb{R}^3}\rho(x)\,dx=c^{-\frac{3}{2}}N$. Consequently, we have the above existence theorem. 
\end{re}
\begin{re}
From \cite{CMP-lieb-Chandrasekhar-theory}, we know that 
for each $0<N<N_*$, the minimizer $\rho^{Ch}$ for  $E^{Ch}(N)$ is unique. Note that  we have the above equivalent existence theorem, thus, for each $0<N^{\frac{2}{3}}<cN_*^{\frac{2}{3}}$, the minimizer $\rho_{c}$ for $E_c(N)$ is also unique. 
\end{re}
 Concerning the nonrelativistic limit $c\to\infty$ of minimizers for $E_c(N)$, we have the following result.
\begin{thm}\label{th1.2}
 Let $\rho_{c}$ be a symmetric decreasing minimizer of $E_{c}(N)$ with fixed $N$ satisfying $0<N^{\frac{2}{3}}<cN_*^{\frac{2}{3}}$. Then, up to  a subsequence $\{\rho_c\}$,  as $c\rightarrow \infty$ we have
\begin{align}\label{eq1.10}
\rho_{c}\to\rho_\infty\ \ \text{strongly  in }L^1\cap L^{\frac{5}{3}}(\mathbb{R}^3), 
\end{align}
where $\rho_\infty$ is a minimizer of $E_\infty(N)$.
\end{thm}

In \cite{Lenz-uniqueness}, the key point to obtain the  nonrelativistic limit of ground
states for the pseudo-relativistic Hartree energy functional  is to use the uniform bound of Lagrange multiplier. However, their method does not work for our model, since the relativistic kinetic energy term $\int_{\mathbb{R}^3}j_{mc}(\rho(x))\,dx$ is nonlinear and  harder to deal with, and  the uniform bound of  Lagrange multiplier $\mu_{c}$ is hard to obtained directly. Fortunately, inspired by the works of  Choi, Seok and Hong in \cite{JFA-convergence-rate}, we find a new way to overcome this problem. We can deduce the uniform lower bound on the operator $H_c:=\sqrt{c^2|p|^2+m^2c^4}-mc^2+\delta$, given by  
$$H_c\ge B|p|,$$ where $B=\min\left\lbrace\frac{2\delta^{1/2}}{\left(2\sqrt{5}m\right)^{1/2}},\frac{c}{2}\right\rbrace$ is a constant. Then by choosing suitable $\delta>0$, we can obtain the uniform bound of $\|\rho_c\|_{L^{4/3}}$, this is the key point to get the uniform  bound of Lagrange multiplier $\mu_{c}$. Finally, combining with $\mu_{c}>0$ we can  obtain the uniform bound of $\|\rho_c\|_{L^{5/3}}$, this is the key point to obtain the nonrelativistic limit \eqref{eq1.10}. Furthermore, we also have the regularity of $\rho_{c}$, in fact $\rho_{c}\in L^1\cap L^{\infty}(\mathbb{R}^3)$.

  For the second part of this paper, we consider the convergence rates of the nonrelativistic limit.  
  
\begin{thm}\label{th1.3}
For any  fixed $N$ satisfying $0<N^{\frac{2}{3}}<cN_*^{\frac{2}{3}}$. Then \\
(i) Let $\rho_{c}(x)$  be a nonnegative minimizer of $E_{c}(N)$ and  $\rho_{\infty}(x)$ be a nonnegative minimizer of $E_{\infty}(N)$,  we have
\begin{align}\label{th1.3-1}
0\leq \int_{\mathbb{R}^3} \rho_{c}(x)^{\frac{5}{3}}\,dx-\int_{\mathbb{R}^3}\rho_{\infty}(x)^{\frac{5}{3}}\,dx=O\left(\frac{1}{c^2}\right).
\end{align}
(ii) Let $R_c$ denotes the radius of the compact support of $\rho_c(x)$ and  $R_\infty$ denotes the  radius of the compact support of $\rho_{\infty}(x)$, we have 
\begin{align}\label{th1.3-2}
 0\leq R_{\infty}-R_{c}=O\left(\frac{1}{c^2}\right).
\end{align}
\end{thm}
\vspace{0.3cm}

From \eqref{th1.3-2}, we know that for $c$ large enough, there exists a constant $K_1>0$ such that
\begin{align}\label{5.13-1}
	B(R_c-\frac{K_1}{c^2})\subset B(R_c)\subset B(R_\infty).
\end{align}
Then we have the following estimates in the corner layer $B(R_\infty)\setminus B(R_c-\frac{K_1}{c^2})$.
\begin{thm} \label{th1.4} 
For $c$ large enough  we have  \\
(i) $\rho_\infty(x)\lesssim \frac{1}{c^3}  $ for all $x\in B(R_\infty)\setminus B(R_c)$. \\
(ii) $\rho_c(x)\lesssim \frac{1}{c^3} $ for all $x\in B(R_c)\setminus B(R_c-\frac{K_1}{c^2})$.
	
\end{thm}

Finally, we give some optimal  uniform bounds  with respect to $N$.
\begin{thm}\label{th1.5}
For any  fixed $N$ satisfying $0<N^{\frac{2}{3}}<cN_*^{\frac{2}{3}}$,
let $\rho_{c}$  be a nonnegative minimizer of $E_{c}(N)$ and let $\rho_{\infty}$ be a nonnegative minimizer of $E_{\infty}(N)$, we have\\
(i) $\|\rho_\infty\|_{L^\infty(\mathbb{R}^3)}\sim N^{2}, \ \  R_{\infty}\sim N^{-\frac{1}{3}}$;\\ 
(ii) $\|\rho_c\|_{L^\infty(\mathbb{R}^3)}\sim N^{2}, \ \ \, R_{c}\sim N^{-\frac{1}{3}}$ for  $c$ large enough.	
	\end{thm}
%2. Since $\rho_c$ and $\rho_\infty$ have compact supports with radii $R_c$ and $R_\infty$, we can not subtract  two corresponding Euler-Lagrange equations directly to obtain the $L^\infty$-convergence rates in Theorem \ref{th1.4}. To overcome this difficult, we split into two parts, in the corner layer $B(R_\infty +\frac{K_1}{c^2})\setminus B(R_\infty -\frac{K_1}{c^2})$ and insider the corner layer $ B(R_\infty -\frac{K_1}{c^2})$. In  the corner layer, we find that $\rho_c$ and $\rho_\infty$ decay fast with the rate $O(\frac{1}{c^3})$, then we can easily obtain the Theorem \ref{th1.4} (i). Inside the corner layer, that is, in $B(R_\infty -\frac{K_1}{c^2})$, we can subtract  two corresponding Euler-Lagrange equations directly, but need to $N$ large enough. This is the most technical part in our paper, we need to get several  estimates with respect to $N$ carefully.  
{\bf Comments:} 

1. By the definition of $j_\infty(\rho)$, \eqref{th1.3-1} is equivalent to $$\int_{\mathbb{R}^3}[j_{\infty}(\rho_{c}(x))-j_{\infty}(\rho_{\infty}(x))]\,dx=O\left(\frac{1}{c^2}\right).$$ For the proof of this, the first difficulty is to deal with the term $D(\rho_c,\rho_c)-D(\rho_\infty,\rho_\infty)$, to overcome  this problem, we establish a novel approach which is depending on  Virial identities to the minimizing problem $E_c(N)$ and $E_{\infty}(N)$ carefully to get some useful energy estimates. In order to prove the convergence rate is optimal, it's vital to get the operator inequality (see Lemma \ref{lem-bound-of-operator}) and the uniformly $L^p$-regularity estimates.

2. The key point to prove \eqref{th1.3-2} is to obtain the estimate $\mu_c-\mu_\infty$. The proof of the upper bound is similar to \eqref{th1.3-1}. However, the lower bound is more complicate, this needs some techniques, the detail for this is presented in section \ref{sec4}.

3. Theorem \ref{th1.4}  shows the decay rates of $\rho_c$ and  $\rho_\infty$ within the corner layer, which can be proved by Newton's theorem.

4. For Theorem \ref{th1.5}, we should mention that Lieb and Yau  \cite{CMP-lieb-Chandrasekhar-theory} proved similar result that the  radius of compact support $R_N\rightarrow \infty$ as $N\rightarrow 0$, but they just did some qualitative analysis,  here we give  concrete best estimates.  
\vspace{0.3cm}

\noindent
\textbf{Organization of the paper:} In section \ref{sec1}, we present the physical background of the equations, status of research and the main results of this paper. In section \ref{sec2}, we give some basic properties. In section \ref{sec3}, we investigate the nonrelativistic limit of minimizer of $E_c(N)$ and  complete the proof of Theorem \ref{th1.2}. In section \ref{sec4}, we consider the convergence rate  and give the proof of Theorem \ref{th1.3} and Theorem \ref{th1.4}. In section \ref{sec5}, we give some optimal estimates with respect to $N$.

\vspace{0.2cm}

\noindent
\textbf{Notation:}\\
- $\rightharpoonup$ denotes weakly converge.\\
- $f*h$  denotes the convolution on $\mathbb{R}^3$.\\
- $\|\cdot\|_{L^s}$ denotes the $L^s(\mathbb{R}^3)$ norm for $s\geq 1$.\\
- $L^p_w(\mathbb{R}^3)$ denotes  weak $L^p$-space, one can see \cite{lieb-Loss-analysis}. It is defined as the space of all measurable functions $f$ such that $$\sup_{\alpha>0}\alpha
\left| \left\{ x:|f(x)|>\alpha \right\} \right|^{1/p}<\infty.$$
- The value of positive constant $C$ is allowed to change from line to line and also in the same formula.\\
- $X\lesssim Y\ (X\gtrsim Y)$ denotes $X\leq CY\ (X\geq CY)$ for some appropriate positive constants $C$ independent of $c$.\\
- $f(c)=O\left(\frac{1}{c^2}\right)$ denotes that there exist two  constants $K_1>0$ and $K_2>0$ such that $\frac{K_1}{c^2}\leq f(c)\leq \frac{K_2}{c^2}.$

\section{Preliminaries}\label{sec2}

Suppose $\rho_{c}$ is a minimizer of $E_c(N)$, $\rho_{\infty}$ is a minimizer of $E_\infty(N)$, by \cite{CMP-lieb-Chandrasekhar-theory} (or \cite{lieb-Loss-analysis}) then $\rho_{c}$ satisfies the following Euler-Lagrange equation
\begin{gather}\label{eq-Euler-equation1}
	\begin{aligned}
		\sqrt{c^2\eta_{c}^2+m^2c^4}-mc^2-\kappa\left(|x|^{-1}*\rho_{c}(x)\right)+\mu_{c}
		\begin{cases}
			=0\quad \text{if}\ \rho_{c}(x)>0,\\
			\geq0\quad \text{if}\ \rho_{c}(x)=0,
		\end{cases}
	\end{aligned}
\end{gather}
or 
\begin{align}\label{eq-Euler-equation1-1}
		\sqrt{c^2\eta_{c}^2+m^2c^4}-mc^2=\left[ \left(\kappa|x|^{-1}*\rho_{c}(x)\right)-\mu_{c}\right]_{+},
\end{align}
where $\eta_c=(6\pi^2\rho_c/q)^{\frac{1}{3}}$, $\mu_{c}$ is a Lagrange multiplier, $[f(x)]_+:=\max\{f(x),0\}$. Similarly,  $\rho_{\infty}$ satisfies
\begin{gather}\label{eq-Euler-equation2}
	\begin{aligned}
		\frac{\eta_{\infty}^2}{2m}-\kappa(|x|^{-1}*\rho_\infty(x))+\mu_{\infty}
		\begin{cases}
			=0, \ \ \text{if } \rho_{\infty}(x)>0,\\
			\ge0, \ \ \text{if } \rho_{\infty}(x)=0,		
		\end{cases}
	\end{aligned}
\end{gather}
or 
\begin{align}\label{eq-Euler-equation2-1}
		\frac{\eta_{\infty}^2}{2m}=\left[ \left(\kappa|x|^{-1}*\rho_{\infty}(x)\right)-\mu_{\infty}\right]_{+},
\end{align}
where $\eta_\infty=(6\pi^2\rho_\infty/q)^{\frac{1}{3}}$, $\mu_{\infty}$ is a Lagrange multiplier.

Firstly, we have
\begin{lem}\label{lem-mul0}
    $\mu_c>0$ and $\mu_{\infty}>0$, furthemore, we have 
\begin{align}\label{eq25}
    \mu_{c}=\frac{\kappa N}{R_{c}}, \quad \mu_{\infty}=\frac{\kappa N}{R_{\infty}},
\end{align}
where $R_c$ and $R_\infty$ denote the radius of the compact support of $\rho_c(x)$ and  $\rho_\infty(x)$.
\end{lem}
\begin{pf}
The similar proof can be found in \cite{{CMP-lieb-Chandrasekhar-theory}}, for the sake of completeness, we present its proof for the situation at hand. Let $V_{\rho_{c}}(x)=\left[\kappa|x|^{-1}*\rho_{c}\right](x)$, since $\rho_{c}$ is symmetric decreasing  function, then by Newton's theorem, we have 
\begin{align}\label{eq2.6}
    \left[|x|^{-1}*\rho_{c}\right](x)=\frac{1}{|x|}\int_{|y|\le|x|}\rho_c(y)\,dy+\int_{|y|\ge|x|}\frac{1}{|y|}\rho_{c}(y)\,dy.
\end{align}
Notice the operator inequality $\sqrt{c^2\eta_{c}^2+m^2c^4}-mc^2\le\frac{\eta_c^2}{2m}$, then, from \eqref{eq-Euler-equation1-1} we have 
\begin{align}\label{eq2.7}
    \left[ \left(\kappa|x|^{-1}*\rho_{c}(x)\right)-\mu_{c}\right]_{+}\le\frac{\eta_c^2}{2m}.
\end{align}
If $\mu_{c}\le0$, combining \eqref{eq2.6} and \eqref{eq2.7}, we have 
\begin{align}
V_{\rho_{c}}\le\frac{\eta_c^2}{2m} =\frac{1}{2m}\cdot\left(\frac{6\pi^2}{q}\right)^{\frac{2}{3}}\rho_{c}^\frac{2}{3}.
\end{align}
From \eqref{eq2.6}, we note that $V_{\rho_{c}}(r)\approx N/r$ as $r\to\infty$, thus $\mu_c>0$, for otherwise for large $r$, $\rho_{c}$ couldn't be in $L^{1}(\mathbb{R}^3)$, which contradicts with $\rho_{c}\in L^{1}(\mathbb{R}^3)$. In the same way, we can deduce $\mu_{\infty}>0$.

Next, we show  that $\mu_{c}=\frac{\kappa N}{R_{c}}$ and $\mu_{\infty}=\frac{\kappa N}{R_{\infty}}$. Review \eqref{eq-Euler-equation1-1}, if $\left[ \left(\kappa|x|^{-1}*\rho_{c}(x)\right)-\mu_{c}\right]_{+}=0$,
then 
$$\mu_c=\kappa|x|^{-1}*\rho_{c}(x).$$
Similarly, by \eqref{eq-Euler-equation2-1}, we have $$\mu_\infty=\kappa|x|^{-1}*\rho_{\infty}(x).$$
Combining with \eqref{eq2.6}, this  implies  \eqref{eq25} holds. This completes the proof of Lemma \ref{lem-mul0}.
\end{pf}

Next, we illustrate the relationship of the operator $H_c=\sqrt{c^2|p|^2+m^2c^4}-mc^2+\delta$ with $|p|$, which plays an important role in  the following sections.
\begin{lem}\label{lem2.1-operator-bounded}
Let $H_c=\sqrt{c^2|p|^2+m^2c^4}-mc^2+\delta$ with $\delta>0$ independent of $c$, then we have 
\begin{align}\label{eq2.1}
	H_c\ge B|p|,
\end{align}
where $B=\min\left\lbrace\frac{2\delta^{1/2}}{\left(2\sqrt{5}m\right)^{1/2}},\frac{c}{2}\right\rbrace$ is a constant which depends  on $m$ and $c$.
\end{lem}
\begin{pf} The similar proof  can be found in \cite{JFA-convergence-rate}, for the sake of completeness, we present the detail.
Factoring out $mc^2$ from the square root in the symbol, we write
\begin{align}\label{eq2.2}
 H_{c}=mc^2(\sqrt{1+|\frac{p}{mc}|^2}-1)+\delta=mc^2f(|\frac{p}{mc}|^2)+\delta,
\end{align}
where $$f(t)=\sqrt{1+t}-1.$$
By the first Taylor expansion, if $0\le t\le4$, then 
\begin{gather}\label{eq2.3}
\begin{aligned}
f(t)=\sqrt{1+t}-1&=f(0)+f'(t_*)t,\text{ for some }t_*\in[0,4]\\
&=\frac{t}{2\sqrt{1+t_*}}\ge\frac{t}{2\sqrt{5}}.
\end{aligned}
\end{gather}
Hence, if $|p|\le2mc$, then
\begin{gather}\label{eq2.4}
	\begin{aligned}
H_{c}&=mc^2(\sqrt{1+|\frac{p}{mc}|^2}-1)+\delta=mc^2f(|\frac{p}{mc}|^2)+\delta\\
	&\ge mc^2\frac{|\frac{p}{mc}|^2}{2\sqrt{5}}+\delta=\frac{|p|^2}{2\sqrt{5}m}+\delta\\
   &\ge\frac{2\delta^{1/2}|p|}{\left(2\sqrt{5}m\right)^{1/2}},
	\end{aligned}
\end{gather}
on the last inequality we use the fact $a^2+b^2\ge 2ab$.

On the other hand, if $|p|\ge 2mc$, then we have 
\begin{gather}\label{eq2.5}
\begin{aligned}
H_{c}&=c|p|\sqrt{1+|\frac{mc}{p}|^2}-mc^2+\delta\\
&\ge c|p|-mc^2+\delta\\
&\ge c|p|-\frac{c|p|}{2}+\delta\geq\frac{c|p|}{2}.
\end{aligned}
\end{gather}
This implies that the lemma holds.
\end{pf}

\section{The nonrelativistic limit }\label{sec3}
In this section, we consider the  nonrelativistic limit. We start with the following lemmas, which estimate two bounds depending on $N$,  which are useful in the section \ref{sec5} below. 
\begin{lem}\label{lem-31boundness}
For fixed $N$ with $0<N^{\frac{2}{3}}<cN_*^{\frac{2}{3}}$, let $\rho_{c}$ be a minimizer of $E_{c}(N)$, then for $c$ large enough, we have 
\begin{align}
	\|\rho_{c}\|_{L^{\frac{4}{3}}(\mathbb{R}^3)}^{\frac{4}{3}}\le \frac{2\sqrt{5}m}{K_{cl}N_*^{\frac{2}{3}}}\, N^{\frac{5}{3}}.
\end{align}
\end{lem}
\begin{pf}
By Lemma \ref{lem2.1-operator-bounded}, we observe that 
\begin{gather}\label{eq27}
	\begin{aligned}
         \frac{3}{4}\left(\frac{6\pi^2}{q} \right)^{\frac{1}{3}}\|\rho_{c}\|_{L^{\frac{4}{3}}}^{\frac{4}{3}}&= \frac{q}{8\pi^3}\int_{\mathbb{R}^3}\int_{|p|<\eta_c}	|p|\,dpdx\\
        &\le\frac{q}{8\pi^3B}\int_{\mathbb{R}^3}\int_{|p|<\eta_c}\left(\sqrt{c^2|p|^2+m^2c^4}-mc^2+\delta\right) \,dpdx,
	\end{aligned}
\end{gather}
where $\eta_c=\left( \frac{6\pi^2\rho_c}{q}\right)^{\frac{1}{3}}$. By \eqref{GNS2} and \eqref{eq27}, we notice that 
\begin{gather}\label{eq2.8}
	\begin{aligned}
\mathcal{E}_c(\rho_c)&=\int_{\mathbb{R}^3}j_{mc}(\rho_{c}(x))\,dx-\frac{\kappa }{2}\int_{\mathbb{R}^3}\int_{\mathbb{R}^3}\frac{\rho_{c}(x)\rho_{c}(y)}{|x-y|}\,dxdy\\
&\ge \left( \frac{3B}{4}\left(\frac{6\pi^2}{q} \right)^{\frac{1}{3}}\|\rho_{c}\|_{L^{\frac{4}{3}}}^{\frac{4}{3}}-\delta N\right) -\kappa D(\rho_{c},\rho_{c})\\
&\ge \frac{3B}{4}\left(\frac{6\pi^2}{q} \right)^{\frac{1}{3}}\|\rho_{c}\|_{L^{\frac{4}{3}}}^{\frac{4}{3}}-\delta N-K_{cl}\left( \frac{N}{N_*}\right)^{\frac{2}{3}}\|\rho_{c}\|_{L^{\frac{4}{3}}}^{\frac{4}{3}}\\
&=K_{cl}(B-\left( \frac{N}{N_*}\right)^{\frac{2}{3}})\|\rho_{c}\|_{L^{\frac{4}{3}}}^{\frac{4}{3}}-\delta N,
	\end{aligned}
\end{gather}
where $K_{cl}=\frac{3}{4}\left({6\pi^2}/{q}\right)^{1/3}$. Note that $B=\min\left\lbrace\frac{2\delta^{1/2}}{\left(2\sqrt{5}m\right)^{1/2}},\frac{c}{2}\right\rbrace$ (see Lemma \ref{lem2.1-operator-bounded}), since $c\to\infty$, take $$\delta=2\sqrt{5}m\left( \frac{N}{N_*}\right)^{\frac{4}{3}},$$
then for $c$ large enough
\begin{align}\label{B1}
B-\left( \frac{N}{N_*}\right)^{\frac{2}{3}}=2\left( \frac{N}{N_*}\right)^{\frac{2}{3}}-\left( \frac{N}{N_*}\right)^{\frac{2}{3}}=\left( \frac{N}{N_*}\right)^{\frac{2}{3}}.
\end{align}
Let $\rho_{\infty}$ be a minimizer of $E_{\infty}(N)$. Using the operator inequality $\sqrt{c^2|p|^2+m^2c^4}-mc^2\le\frac{|p|^2}{2m}$, the definition of $\rho_{c}$ and the fact that $E_\infty(N)<0$ (see Appendix A), we have 
\begin{align}\label{eq:enengy ralation}
	\mathcal{E}_c(\rho_{c})\le\mathcal{E}_{c}(\rho_\infty)\le\mathcal{E}_{\infty}(\rho_\infty)=E_{\infty}(N)<0.
\end{align}
Combining \eqref{eq2.8}, \eqref{B1} and  \eqref{eq:enengy ralation}, we deduce that
\begin{align*}
K_{cl}\left( \frac{N}{N_*}\right)^{\frac{2}{3}}\|\rho_{c}\|_{L^{\frac{4}{3}}}^{\frac{4}{3}}\le\mathcal{E}_c(\rho_{c})+\delta N\le E_{\infty}(N)+\delta N\leq \delta N.	
\end{align*}
Thus, we have 
\begin{align}\label{eq2.11}
\|\rho_{c}\|_{L^{\frac{4}{3}}}^{\frac{4}{3}}\le \frac{2\sqrt{5}m}{K_{cl}N_*^{\frac{2}{3}}} N^{\frac{5}{3}}.
\end{align}
This completes the proof of the lemma.
\end{pf}

\begin{lem}\label{lem-32bounded}
For fixed $N$ with $0<N^{\frac{2}{3}}<cN_*^{\frac{2}{3}}$, let $\rho_{c}$ be a minimizer of $E_{c}(N)$, then there exists a constant $M>0$ such that for $c$ large enough 
\begin{align}\label{rho-c-53-N}
	\|\rho_{c}\|_{L^{\frac{5}{3}}(\mathbb{R}^3)}\le M N^{\frac{7}{5}},
\end{align}
where the constant $M$ is independent of $c$.
\end{lem}
\begin{pf}
It follows from \eqref{eq-Euler-equation1} that
\begin{gather}
	\begin{aligned}
		c^2&A_0\|\rho_{c}\|_{L^{\frac{5}{3}}}^{\frac{5}{3}}+m^2c^4N\\
      =&\left\langle \sqrt{c^2\eta_{c}^2+m^2c^4}\,\rho_{c}^{\frac{1}{2}},\sqrt{c^2\eta_{c}^2+m^2c^4}\,\rho_{c}^{\frac{1}{2}}\right\rangle\\
		=&\left\langle(-\mu_{c}+\kappa(|x|^{-1}*\rho_{c})+mc^2)\,\rho_{c}^{\frac{1}{2}},(-\mu_{c}+\kappa(|x|^{-1}*\rho_{c}+mc^2))\,\rho_{c}^{\frac{1}{2}}\right\rangle\\
		=&\mu_{c}^2N-2\mu_{c}\kappa\left\langle\rho_{c}^{\frac{1}{2}},(|x|^{-1}*\rho_{c})\rho_{c}^{\frac{1}{2}}\right\rangle+2mc^2\kappa\left\langle\rho_{c}^{\frac{1}{2}},(|x|^{-1}*\rho_{c})\rho_{c}^{\frac{1}{2}}\right\rangle\\
	 &+\kappa^2\left\langle(|x|^{-1}*\rho_{c})\rho_{c}^{\frac{1}{2}},(|x|^{-1}*\rho_{c})\rho_{c}^{\frac{1}{2}} \right\rangle	-2\mu_{c}mc^2N+m^2c^4N.
	\end{aligned}
\end{gather}
where $A_0=(6\pi^2/q)^{\frac{2}{3}}$. Since $\mu_{c}>0$, then we have
\begin{gather}\label{eq3.9-1}
	\begin{aligned}
	c^2A_0\|\rho_{c}\|_{L^{\frac{5}{3}}}^{\frac{5}{3}}\le&\mu_{c}^2N+2mc^2\kappa\left\langle\rho_{c}^{\frac{1}{2}},(|x|^{-1}*\rho_{c})\rho_{c}^{\frac{1}{2}}\right\rangle\\
&+\kappa^2\left\langle(|x|^{-1}*\rho_{c})\rho_{c}^{\frac{1}{2}},(|x|^{-1}*\rho_{c})\rho_{c}^{\frac{1}{2}} \right\rangle.
	\end{aligned}
\end{gather}
By the Hardy-Littlewood-Sobolev inequality and the Interpolation inequality, we obtain
\begin{gather}\label{eq218}
	\begin{aligned}
		\left\langle\rho_{c}^{\frac{1}{2}},(|x|^{-1}*\rho_{c})\rho_{c}^{\frac{1}{2}}\right\rangle\le\,C_1\|\rho_{c}\|_{L^{\frac{6}{5}}}^2\le\,C_1\|\rho_{c}\|_{L^1}^{\frac{7}{6}}\|\rho_{c}\|_{L^{\frac{5}{3}}}^{\frac{5}{6}}=C_1 N^{\frac{7}{6}}\|\rho_{c}\|_{L^{\frac{5}{3}}}^{\frac{5}{6}}.
	\end{aligned}
\end{gather}
Using H$\ddot{\text{o}}$lder's inequality, we have
\begin{gather}
	\begin{aligned}
	\left\langle(|x|^{-1}*\rho_{c})\,\rho_{c}^{\frac{1}{2}},(|x|^{-1}*\rho_{c})\,\rho_{c}^{\frac{1}{2}}\right\rangle\le\|(|x|^{-1}*\rho_{c})^2\|_{L^{\frac{5}{2}}}\|\rho_{c}\|_{L^{\frac{5}{3}}}.
	\end{aligned}
\end{gather}
Notice that $|x|^{-1}\in L_\omega^3$, by the weak Young inequality, we have
\begin{gather}\label{eq220}
	\begin{aligned}
		\|(|x|^{-1}*\rho_{c})^2\|_{L^{\frac{5}{2}}}&=\||x|^{-1}*\rho_{c}\|_{L^5}^2\\
		&\le\||x|^{-1}\|_{L_{\omega}^3}^2\|\rho_{c}\|_{L^{\frac{15}{13}}}^2\\
		&\le\,C_2\|\rho_{c}\|_{L^1}^{\frac{4}{3}}\|\rho_{c}\|_{L^{\frac{5}{3}}}^{\frac{2}{3}}=C_2 N^{\frac{4}{3}}\|\rho_{c}\|_{L^{\frac{5}{3}}}^{\frac{2}{3}}.
	\end{aligned}
\end{gather}
On the last inequality we use the interpolation inequality.  Thus 
\begin{gather}\label{eq313-1}
	\begin{aligned}
	\left\langle(|x|^{-1}*\rho_{c})\,\rho_{c}^{\frac{1}{2}},(|x|^{-1}*\rho_{c})\,\rho_{c}^{\frac{1}{2}}\right\rangle\le C_2 N^{\frac{4}{3}}\|\rho_{c}\|_{L^{\frac{5}{3}}}^{\frac{5}{3}}.
	\end{aligned}
\end{gather}

 On the other hand, by \eqref{eq-Euler-equation1}, we have 
 \begin{align}\label{eq221}
 	-\mu_{c} N=\int_{\mathbb{R}^3}(\sqrt{c^2\eta_{c}^2+m^2c^4}-mc^2)\rho_{c}(x)\,dx-2\kappa\,D(\rho_{c},\rho_{c}).
 \end{align}
Combining $\mu_{c}>0$ with $\sqrt{c^2\eta_{c}^2+m^2c^4}-mc^2>0$, then, by \eqref{GNS2}, we obtain
\begin{align}
\mu_{c}N\le2\kappa\,D(\rho_{c},\rho_{c})\le2K_{cl}\left(\frac{N}{N_*}\right)^{\frac{2}{3}}\|\rho_{c}\|_{L^{\frac{4}{3}}}^{\frac{4}{3}}.
\end{align}
By Lemma \ref{lem-31boundness}, it follows that
\begin{align}\label{eq223}
	\mu_{c}\le 2\delta =\frac{4\sqrt{5}m}{N_*^{\frac{4}{3}}}\,N^{\frac{4}{3}}.
\end{align}
Combining \eqref{eq3.9-1}, \eqref{eq218}, \eqref{eq313-1} and \eqref{eq223}, we deduce that
\begin{align}\label{eq224}
c^2A_0\|\rho_{c}\|_{L^{\frac{5}{3}}}^{\frac{5}{3}}\leq \frac{80 m^2}{N_*^{\frac{8}{3}}}N^{\frac{11}{3}}+2mc^2\kappa  C_1 N^{\frac{7}{6}}\|\rho_{c}\|_{L^{\frac{5}{3}}}^{\frac{5}{6}}+\kappa^{2}C_2 N^{\frac{4}{3}}\|\rho_{c}\|_{L^{\frac{5}{3}}}^{\frac{5}{3}},
\end{align}
Since $N^{\frac{2}{3}}<cN^{\frac{2}{3}}_*$, then $\frac{N^{\frac{4}{3}}}{N^{\frac{4}{3}}_*c^2}<1$, it follows that for $c$ large enough
\begin{align*}
	\|\rho_{c}\|_{L^{\frac{5}{3}}}^{\frac{5}{3}}&\leq \frac{80m^2}{N_*^{\frac{8}{3}}c^2A_0}N^{\frac{11}{3}}+\frac{2m\kappa C_1}{A_0}N^{\frac{7}{6}}\|\rho_{c}\|_{L^{\frac{5}{3}}}^{\frac{5}{6}}+\frac{\kappa^{2}C_2 }{c^2A_0}N^{\frac{4}{3}}\|\rho_{c}\|_{L^{\frac{5}{3}}}^{\frac{5}{3}}\\
	&\leq  \frac{80m^2}{N_*^{\frac{4}{3}}A_0}N^{\frac{7}{3}}+\frac{2m\kappa C_1}{A_0}N^{\frac{7}{6}}\|\rho_{c}\|_{L^{\frac{5}{3}}}^{\frac{5}{6}}+\frac{\kappa^{2}C_2 }{c^2A_0}N^{\frac{4}{3}}\|\rho_{c}\|_{L^{\frac{5}{3}}}^{\frac{5}{3}}\\
	&\leq \frac{80m^2}{N_*^{\frac{4}{3}}A_0}N^{\frac{7}{3}} +\frac{2m\kappa C_1}{A_0}N^{\frac{7}{6}}\|\rho_{c}\|_{L^{\frac{5}{3}}}^{\frac{5}{6}}+\frac{1 }{2} \|\rho_{c}\|_{L^{\frac{5}{3}}}^{\frac{5}{3}}   \ \ \  (\text{$c$ large enough}).
\end{align*}
Thus  we have
$$\frac{1}{2}\|\rho_{c}\|_{L^{\frac{5}{3}}}^{\frac{5}{3}}\leq \frac{80m^2}{N_*^{\frac{4}{3}}A_0} N^{\frac{7}{3}}+\frac{2m\kappa C_1}{A_0} N^{\frac{7}{6}}\|\rho_{c}\|_{L^{\frac{5}{3}}}^{\frac{5}{6}}.$$
This implies that 
$$\|\rho_{c}\|_{L^{\frac{5}{3}}}^{\frac{5}{3}} \leq M N^{\frac{7}{3}},$$
 where  $M>0$ is a constant  independent of $N$ and $c$. Thus \eqref{rho-c-53-N} holds, we complete the proof of Lemma \ref{lem-32bounded}.
\end{pf}

{\bf{The end proof of Theorem \ref{th1.2}:}}
Firstly, we claim that $\{\rho_{c}\}$ is a minimizing sequence of $E_\infty(N)$. Note that
\begin{gather}\label{eq318}
\begin{aligned}E_\infty(N)-E_c(N)&=\mathcal{E}_\infty(\rho_{\infty})-\mathcal{E}_{c}(\rho_{c})\geq \mathcal{E}_\infty(\rho_{\infty})-\mathcal{E}_{c}(\rho_{\infty})\\
&=\frac{q}{8\pi^3}\int_{\mathbb{R}^3}\int_{|p|<\eta_\infty}\left[\frac{|p|^2}{2m}-(\sqrt{c^{2}|p|^2 +m^2c^{4}}-mc^2)\right]\,dpdx\geq 0.
\end{aligned}
\end{gather}
On the other hand, by the operator inequality $\sqrt{c^2|p|^2 +m^2c^{4}}\le\frac{|p|^2}{2m}+mc^2$, we have
\begin{gather}
	\begin{aligned} E_\infty(N)-E_{c}(N)&\le\mathcal{E}_\infty(\rho_{c})-\mathcal{E}_{c}(\rho_{c})  \\
	&=\frac{q}{8\pi^3}\int_{\mathbb{R}^3}\int_{|p|<\eta_c}\left[\frac{|p|^2}{2m}-(\sqrt{c^{2}|p|^2 +m^2c^{4}}-mc^2)\right]\,dpdx\\
	&\leq \frac{q}{8\pi^3}\int_{\mathbb{R}^3}\int_{|p|<\eta_c}\left(\frac{|p|^2}{2m}-\frac{|p|^2}{\sqrt{\frac{|p|^2}{c^2}+m^2}+m}\right)\,dpdx\\
	&\le2h\int_{\mathbb{R}^3}\rho_{c}(x)^{\frac{5}{3}}\,dx,
	\end{aligned}	
\end{gather}
where $h=\frac{3}{10m}(6\pi^2/q )^{2/3}$. From Lemma \ref{lem-32bounded}, we know that  $\int_{\mathbb{R}^3}\rho_{c}(x)^{\frac{5}{3}}\,dx$ is uniformly bounded. Notice  that
 \begin{align}\label{eq-oprator-convergence}
 \sqrt{c^2|p|^2  +m^2c^4}-mc^2-\frac{|p|^2}{2m}\to0 \quad\text{for every }p\in \mathbb{R}^3\text{ as }c\to\infty.
 \end{align}
 Combining  the dominated convergence theorem and \eqref{eq318}-\eqref{eq-oprator-convergence}, we deduce that
 \begin{align}\label{eq321}
 E_{c}(N)\to E_\infty(N)\qquad\text{ and }\quad\mathcal{E}_\infty(\rho_c)\to E_\infty(N).
 \end{align}
 Hence, $\left\{{\rho_{c}}\right\}$ is a minimizing sequence of $E_\infty(N)$. Combining with the existence of minimizer of  $\mathcal{E}_\infty(\rho)$ (see Appendix A), we deduce that \eqref{eq1.10}  holds.{{\hfill $\square$}\\}
 
 \section{The convergence rate }\label{sec4}
Let $\rho_{c}$ be a minimizer of $\mathcal{E}_{c}(\rho)$, then $\rho_{c}$ satisfies the Euler-Lagrange equation 
 \begin{gather}
 	\begin{aligned}
 \sqrt{c^2\eta_{c}^2+m^2c^4}-mc^2-\kappa\left(|x|^{-1}*\rho_{c}(x)\right)+\mu_{c}
 		\begin{cases}
 			=0\quad \text{if}\ \rho_{c}(x)>0,\\
 			\geq0\quad \text{if}\ \rho_{c}(x)=0,
 		\end{cases}
 	\end{aligned}
 \end{gather}
 where $\eta_c=(6\pi^2\rho_c/q)^{\frac{1}{3}}$, $\mu_{c}$ is a Lagrange multiplier. In fact, we have
 \begin{align}\label{eq-mu_c}
 -\mu_{c} N=\int_{\mathbb{R}^3}(\sqrt{c^2\eta_{c}^2+m^2c^4}-mc^2)\rho_{c}(x)\,dx-2\kappa\,D(\rho_{c},\rho_{c}).
 \end{align}
 Moreover, we have the Virial identity $\left.\partial_\lambda\mathcal{E}_{c}(\lambda^3\rho_{c}(\lambda x))\right|_{\lambda=1}=0$, i.e.,
 \begin{gather}\label{eq-Pohozaev-type-identity1}
 	\begin{aligned}
 		\int_{\mathbb{R}^3}j_{mc}(\rho_{c}(x))\,dx-m^2c^4\int_{\mathbb{R}^3}\bar{j}_{mc}(\rho_{c}(x))\,dx-\kappa D(\rho_{c},\rho_{c})+mc^2N=0,
 	\end{aligned}
 \end{gather}
 where
 \begin{gather}\label{jmcbarj}
 	\begin{aligned} 
 		\bar{j}_{mc}(\rho_c)&=\frac{q}{8\pi^3}\int_{|p|<(6\pi^2\rho_c/q)^{1/3}}\frac{1}{\sqrt{c^2|p|^2+m^2c^4}}\,dp\\
 		&=\frac{q}{4\pi^2c}\left[\eta_c\sqrt{\eta_c^2+m^2c^2}-m^2c^2\ln\left(\frac{\eta_c+\sqrt{\eta_c^2+m^2c^2}}{mc}\right)\right],\eta_c=\left(\frac{6\pi^2\rho_c}{q}\right)^{1/3}.
 	\end{aligned}
 \end{gather}
 From \eqref{eq-Pohozaev-type-identity1}, we have
 \begin{align}\label{eq-Drhoc}
 	-\kappa D(\rho_{c},\rho_{c})=2\mathcal{E}_c(\rho_c)-\int_{\mathbb{R}^3}j_{mc}(\rho_{c}(x))\,dx-m^2c^4\int_{\mathbb{R}^3}\bar{j}_{mc}(\rho_{c}(x))\,dx+mc^2N.
 \end{align}
 
 Let $\rho_\infty$ be a minimizer of $\mathcal{E}_\infty(\rho)$, then it satisfies  Euler-Lagrange equation 
 \begin{gather}
 	\begin{aligned}
 		\frac{\eta_{\infty}^2}{2m}-\kappa(|x|^{-1}*\rho_\infty(x))+\mu_{\infty}
 		\begin{cases}
 			=0 \ \ \text{if } \rho_{\infty}(x)>0,\\
 			\ge0 \ \ \text{if } \rho_{\infty}(x)=0,		
 		\end{cases}
 	\end{aligned}
 \end{gather}
 where $\eta_{\infty}=(6\pi^2\rho_{\infty}/q)^{\frac{1}{3}}$ and $\mu_{\infty}$ is a Lagrange multiplier. In fact, we also have
 \begin{gather}\label{eq-mu-infty}
 \begin{aligned}
 -\mu_{\infty}N&=\frac{1}{2m}\int_{\mathbb{R}^3}\eta_{\infty}^2\rho_\infty\,dx-2\kappa D(\rho_\infty,\rho_\infty)\\
&=\frac{1}{2m}\big(\frac{6\pi^2}{q}\big)^{\frac{2}{3}}\int_{\mathbb{R}^3}\rho_\infty(x)^{\frac{5}{3}}\,dx-2\kappa D(\rho_\infty,\rho_\infty).
 \end{aligned}
 \end{gather}
 Moreover, we have the Virial identity 
 $\partial_\lambda \mathcal{E}_{\infty} (\lambda ^3\rho_{\infty}(\lambda x))|_{\lambda =1}=0$, i.e.,
 \begin{gather}\label{eq-pohozaev-identity2}
 	\begin{aligned}
 2\int_{\mathbb{R}^3}j_{\infty}(\rho_{\infty}(x))\,dx-\kappa D(\rho_\infty,\rho_\infty)=\frac{3}{5m}\big(\frac{6\pi^2}{q}\big)^{\frac{2}{3}}\int_{\mathbb{R}^3}\rho_\infty(x)^{\frac{5}{3}}\,dx-\kappa D(\rho_\infty,\rho_\infty)=0,
 	\end{aligned}
 \end{gather}
 where $j_{\infty}(\rho_{\infty}(x))$ is defined in \eqref{j-infty-rho}.
 By \eqref{eq-pohozaev-identity2}, we have 
 \begin{align}\label{eq-Drhoinfty}
 	-\kappa D(\rho_\infty,\rho_\infty)=2\mathcal{E}_{\infty} (\rho_\infty).
 \end{align}
\subsection{Energy estimates}
Firstly, we have the following uniform $L^\infty$ estimate.
\begin{lem}\label{lem2.3-rhoc-regularity}
 There exists a constant $M_1>0$ is dependent of $N$ and independent of $c$ such that $$||\rho_{c}||_{L^\infty(\mathbb{R}^3)}\leq M_1.$$
\end{lem}
\begin{pf}
Note that by direct calculation we have 
\begin{gather}\label{eq2.15}
 \begin{aligned}
 |x|^{-1}*\rho_{c}(x)&=\int_{|x-y|\le1}\frac{1}{|x-y|}\rho_{c}(y)\,dy+\int_{|x-y|>1}\frac{1}{|x-y|}\rho_{c}(y)\,dy\\
  &\le\int_{|z|\le1}\frac{1}{|z|}\rho_{c}(z+x)\,dz+\int_{\mathbb{R}^3}\rho_{c}(y)\,dy\\
  &\le\left(\int_{|z|\le1}\frac{1}{|z|^\frac{5}{2}}\,dz\right)^{\frac{2}{5}}\left( \int_{|z|\le1}\rho_{c}^{\frac{5}{3}}(z+x)\,dz\right)^{\frac{3}{5}}+N\\
  &\lesssim\|\rho_{c}\|_{L^{\frac{5}{3}}}+N.
  \end{aligned}
	\end{gather}
By Lemma \ref{lem-32bounded}, we know that $\|\rho_{c}\|_{L^{\frac{5}{3}}}$ is uniformly bounded, therefore, we deduce that there exists a constant $M_2>0$ which is independent of $c$ such that
\begin{align}\label{eq226}
	\||x|^{-1}*\rho_{c}(x)\|_{L^{\infty}}\le M_2.
\end{align}

Next, we rewrite \eqref{eq-Euler-equation1-1} as 
\begin{align}
		\sqrt{c^2\eta_{c}^2+m^2c^4}-mc^2+1=\left[ \left(\kappa|x|^{-1}*\rho_{c}(x)\right)-\mu_{c}\right]_{+}+1.
\end{align}
Let $\delta=1$ in Lemma \ref{lem2.1-operator-bounded}, we have $\sqrt{c^2\eta_{c}^2+m^2c^4}-mc^2+1\ge B|\eta_{c}|$, thus we deduce that 
\begin{align}\label{eq4.4}
\left[\left(\kappa|x|^{-1}*\rho_{c}(x)\right)-\mu_{c}\right]_{+}+1\ge B|\eta_{c}|,
\end{align}
where $B=\min\left\lbrace\frac{2}{\left(2\sqrt{5}m\right)^{1/2}},\frac{c}{2}\right\rbrace$, since $c>1$ sufficiently large, then $B=\frac{2}{\left(2\sqrt{5}m\right)^{1/2}}$
is independent of $c$. It follows \eqref{eq223} and \eqref{eq226} that
\begin{align}\label{eq4.5}
    \|\left[\left(\kappa|x|^{-1}*\rho_{c}(x)\right)-\mu_{c}\right]_{+}+1\|_{L^{\infty}}\le M_3,
\end{align}
with $M_3>0$ is independent of $c$. 
Consequently, combining \eqref{eq4.4} with \eqref{eq4.5} we deduce 
\begin{align}
\|\eta_{c}\|_{L^{\infty}}=\left({6\pi^2}/{q}\right)^{\frac{1}{3}}\|\rho_{c}\|_{L^{\infty}}^{\frac{1}{3}}\le M_4,
\end{align}
where $M_4>0$ is a constant which is independent of $c$. This implies the Lemma \ref{lem2.3-rhoc-regularity} holds.
\end{pf}
\begin{lem}\label{lem-regulatity-rhoinfty}
There exists a constant $G>0$ independent of $c$ such that 
\begin{align*}
\|\rho_{\infty}\|_{L^{\infty}(\mathbb{R}^3)}\le G.
\end{align*}
\end{lem}
\begin{pf}
Firstly, we show $\|\rho_{\infty}\|_{L^{\frac{5}{3}}(\mathbb{R}^3)}$ is uniformly bounded. Multiplying \eqref{eq-Euler-equation2-1} with $\rho_{\infty}$, since $\mu_{\infty}>0$, then
\begin{gather}\label{eq416}
\begin{aligned}
\frac{1}{2m}\left(\frac{6\pi^2}{q}\right)^{\frac{2}{3}}\|\rho_{\infty}\|_{L^{\frac{5}{3}}(\mathbb{R}^3)}^{\frac{5}{3}}&=\kappa\int_{\mathbb{R}^3}\int_{\mathbb{R}^3}\frac{\rho_{\infty}(x)\rho_{\infty}(y)}{|x-y|}\,dxdy-\mu_{\infty}N\\
&\le\kappa\int_{\mathbb{R}^3}\int_{\mathbb{R}^3}\frac{\rho_{\infty}(x)\rho_{\infty}(y)}{|x-y|}\,dxdy\\
&\le\kappa\|\rho_{\infty}\|_{L^{\frac{6}{5}}}^{2}\\
&\le\kappa\|\rho_{\infty}\|_{L^1}^{\frac{7}{6}}\|\rho_{\infty}\|_{L^{\frac{5}{3}}}^{\frac{5}{6}}=\kappa N^{\frac{7}{6}}\|\rho_{\infty}\|_{L^{\frac{5}{3}}}^{\frac{5}{6}},
\end{aligned}
\end{gather}
where we use the Hardy-Littlewood-Sobolev inequality and Interpolation inequality. Then, 
\begin{align}\label{eq417}
\|\rho_{\infty}\|_{L^{\frac{5}{3}}(\mathbb{R}^3)}\le CN^{\frac{7}{5}}.
\end{align}

Next, we claim that $\rho_{\infty}\in L^{\infty}(\mathbb{R}^3)$. The similar to \eqref{eq2.15},
we have\begin{gather}\label{eq418}
 \begin{aligned}
 |x|^{-1}*\rho_{\infty}(x)&=\int_{|x-y|\le1}\frac{1}{|x-y|}\rho_{\infty}(y)\,dy+\int_{|x-y|>1}\frac{1}{|x-y|}\rho_{\infty}(y)\,dy\\
&\le\int_{|z|\le1}\frac{1}{|z|}\rho_{\infty}(z+x)\,dz+\int_{\mathbb{R}^3}\rho_{\infty}(y)\,dy\\
&\le\left(\int_{|z|\le1}\frac{1}{|z|^\frac{5}{2}}\,dz\right)^{\frac{2}{5}}\left( \int_{|z|\le1}\rho_{\infty}^{\frac{5}{3}}(z+x)\,dz\right)^{\frac{3}{5}}+N\\
&\lesssim\|\rho_{\infty}\|_{L^{\frac{5}{3}}}+N.
\end{aligned}
\end{gather}
It follows \eqref{eq417} that
\begin{align}
\|\,|x|^{-1}*\rho_{\infty}(x)\|_{L^\infty}\le C_2,
\end{align}
with $C_2>0$ independent of $c$. Finally, by \eqref{eq-Euler-equation2-1}, we have
\begin{align}
\frac{1}{2m}\|\eta_{\infty}^2\|_{L^\infty}=\frac{1}{2m}\left(\frac{6\pi^2}{q}\right)^{\frac{2}{3}}\|\rho_{\infty}\|_{L^\infty}^\frac{2}{3}\le\|\,|x|^{-1}*\rho_{\infty}(x)\|_{L^\infty}\le\,C_2.
\end{align}
This implies the lemma holds.
\end{pf}
\begin{lem}\label{lem-bound-of-operator}
For $c>1$, we have
\begin{align*}
\frac{|p|^4}{8m(|p|^2+m^2)}\cdot\frac{1}{c^2}	\leq \frac{|p|^2}{2m}-\left(\sqrt{c^2|p|^2+m^2c^4}-mc^2\right)\leq\frac{|p|^4}{8m^3}\cdot\frac{1}{c^2}
\end{align*}
\end{lem}
\begin{pf}
Observe that
\begin{align*}
&\frac{|p|^2}{2m}-\left(\sqrt{c^2|p|^2+m^2c^4}-mc^2\right) 
=\frac{|p|^2}{2m}\left(1-\frac{2m}{\sqrt{\frac{|p|^2}{c^2}+m^2}+m}\right)\\
=&\frac{|p|^2}{2m}\left( \frac{\sqrt{\frac{|p|^2}{c^2}+m^2}-m}{\sqrt{\frac{|p|^2}{c^2}+m^2}+m}\right)=\frac{|p|^2}{2m}\left( \frac{\frac{|p|^2}{c^2}}{(\sqrt{\frac{|p|^2}{c^2}+m^2}+m)^2}\right)\\
 =&\frac{|p|^4}{2m}\left( \frac{1}{(\sqrt{\frac{|p|^2}{c^2}+m^2}+m)^2}\right) \cdot\frac{1}{c^2}.
\end{align*}
Since $c>1$, then 
$$2m\leq \sqrt{\frac{|p|^2}{c^2}+m^2}+m\leq 2\sqrt{|p|^2+m^2}.$$
This implies that this lemma holds.
\end{pf}

From \eqref{eq321}, we know that $E_c(N)\to E_\infty(N)$ as $c\to\infty$, the following lemma will show that the above energy convergence has the optimal convergence order with $\frac{1}{c^2}$.

\begin{lem}\label{lem-energy-estimate}
 let $\rho_{c}$  be a nonnegative minimizer of $E_{c}(N)$ and let $\rho_{\infty}$ be a nonnegative minimizer of $E_{\infty}(N)$, then 
 $$E_{\infty}(N)-E_{c}(N)=\mathcal{E}_{\infty}(\rho_{\infty})-\mathcal{E}_{c}(\rho_{c})\sim\frac{1}{c^2}.$$
\end{lem}
\begin{pf}
We start with the upper bound. By the definition of $E_\infty(N)$ and $E_{c}(N)$, we have
\begin{gather}\label{eq421}
 \begin{aligned}
 E_\infty(N)-E_{c}(N)&\le\mathcal{E}_\infty(\rho_{c})-\mathcal{E}_{c}(\rho_{c})  \\
 &=\frac{q}{8\pi^3}\int_{\mathbb{R}^3}\int_{|p|<\eta_c}\left[\frac{|p|^2}{2m}-(\sqrt{c^{2}|p|^2 +m^2c^{4}}-mc^2)\right]\,dpdx.
 \end{aligned}
 \end{gather}
It follows from Lemma \ref{lem-bound-of-operator} that
 \begin{gather}\label{eq422}
  \begin{aligned}
     E_\infty(N)-E_{c}(N)&\le\frac{q}{8\pi^3}.\frac{1}{8m^3c^2}\int_{\mathbb{R}^3}\int_{|p|<\eta_c}|p|^4\,dpdx\\
     &=\frac{A}{c^2}\int_{\mathbb{R}^3}\rho_c(x)^\frac{7}{3}\,dx,
  \end{aligned}
\end{gather}
 where $A=\frac{3}{56m^3}(6\pi^2/q)^{\frac{4}{3}}$. According to Lemma \ref{lem2.3-rhoc-regularity}, we get
 \begin{align}\label{eq423}
 E_\infty(N)-E_{c}(N)\lesssim\frac{1}{c^2}.
 \end{align}

 Next, we show the lower bound. According to \eqref{eq318} and Lemma \ref{lem-bound-of-operator}, we get 
 \begin{gather}\label{eq424}
\begin{aligned}
E_\infty(N)-E_c(N)&=\mathcal{E}_\infty(\rho_{\infty})-\mathcal{E}_{c}(\rho_{c})\\
&\geq\mathcal{E}_\infty(\rho_{\infty})-\mathcal{E}_{c}(\rho_{\infty})\\
&=\frac{q}{8\pi^3}\int_{\mathbb{R}^3}\int_{|p|<\eta_\infty}\left[\frac{|p|^2}{2m}-(\sqrt{c^{2}|p|^2 +m^2c^{4}}-mc^2)\right]\,dpdx\\
&\ge\frac{1}{c^2}\cdot\frac{q}{8\pi^3}\int_{\mathbb{R}^3}\int_{|p|<\eta_\infty}\frac{|p|^4}{8m(|p|^2+m^2)}\,dpdx\\
&\ge\frac{1}{c^2}\cdot\frac{q}{8\pi^3}\int_{\mathbb{R}^3}\int_{|p|<\eta_\infty}\frac{|p|^4}{8m(\eta_{\infty}^2+m^2)}\,dpdx\\
&=\frac{K}{c^2}\int_{\mathbb{R}^3}
\frac{\rho_\infty(x)^\frac{7}{3}}{8m(\eta_{\infty}^2+m^2)}\,dx,
\end{aligned}  
\end{gather}
where $K=\frac{3}{7}(6\pi^2/q)^{\frac{4}{3}}$ and $\eta_{\infty}=(6\pi^2\rho_\infty/q)^{\frac{1}{3}}$.  Then, Lemma \ref{lem-regulatity-rhoinfty} leads to
\begin{align*}
E_\infty(N)-E_c(N)\gtrsim\frac{1}{c^2}.
\end{align*}
This completes the proof of the lemma.
\end{pf}
\subsection{The optimal convergence rate of minimizers}  
In this subsection, we study the convergence rate of minimizers and complete the proof of Theorem \ref{th1.3} (i).

Noticing that for any function $0\le\rho_{c}\in L^{1}\cap L^{\frac{5}{3}}(\mathbb{R}^3)$ with $\int_{\mathbb{R}^3} \rho_{c}(x)\,dx=N$, we have 
$\mathcal{E}_\infty((\rho_{\infty}+s\rho_{c})/(1+s))(s\ge0)$ attains minimum at $s=0$. By direct computation, we find 
\begin{gather}
\begin{aligned}
&\left.\frac{d}{ds}\mathcal{E}_\infty\left(\frac{\rho_{\infty}+s\rho_{c}}{1+s} \right)\right|_{s=0}\\
&=\frac{1}{2m}\cdot\big(\frac{6\pi^2}{q}\big)^{\frac{2}{3}}\int_{\mathbb{R}^3}\rho_\infty(x)^\frac{2}{3}(\rho_c(x)-\rho_{\infty}(x))\,dx-2\kappa D(\rho_{\infty},\rho_{c}-\rho_{\infty})\ge0.
\end{aligned}
\end{gather}
A simple calculation shows
\begin{gather}\label{eq4.26}
\begin{aligned}
\mathcal{E}_\infty&(\rho_{c})-\mathcal{E}_\infty(\rho_{\infty})\\
=&\frac{3}{5}\left.\frac{d}{ds}\mathcal{E}_\infty\left(\frac{\rho_{\infty}+s\rho_{c}}{1+s} \right)\right|_{s=0}\\
&+\frac{3}{10m}\cdot\big(\frac{6\pi^2}{q}\big)^{\frac{2}{3}}\int_{\mathbb{R}^3}\rho_{c}(x)^\frac{5}{3}\,dx-\frac{3}{10m}\cdot\big(\frac{6\pi^2}{q}\big)^{\frac{2}{3}}\int_{\mathbb{R}^3}\rho_{\infty}(x)^\frac{2}{3}\rho_{c}(x)\,dx\\
&+D(\rho_{c},\rho_{c})-D(\rho_{\infty},\rho_{\infty})+\frac{6}{5}D(\rho_{\infty},\rho_{c}-\rho_{\infty})\\
\ge&\frac{3}{10m}\cdot\big(\frac{6\pi^2}{q}\big)^{\frac{2}{3}}\int_{\mathbb{R}^3}\rho_{c}(x)^\frac{5}{3}\,dx-\frac{3}{10m}\cdot\big(\frac{6\pi^2}{q}\big)^{\frac{2}{3}}\int_{\mathbb{R}^3}\rho_{\infty}(x)^\frac{2}{3}\rho_{c}(x)\,dx\\
&+D(\rho_{c},\rho_{c})-D(\rho_{\infty},\rho_{\infty})+\frac{6}{5}D(\rho_{\infty},\rho_{c}-\rho_{\infty})\\
=&\frac{3}{10m}\cdot\big(\frac{6\pi^2}{q}\big)^{\frac{2}{3}}\int_{\mathbb{R}^3}\rho_{c}(x)^\frac{5}{3}\,dx-\frac{3}{10m}\cdot\big(\frac{6\pi^2}{q}\big)^{\frac{2}{3}}\int_{\mathbb{R}^3}\rho_{\infty}(x)^\frac{2}{3}\rho_{c}(x)\,dx\\
&+D(\rho_{c},\rho_{c})-D(\rho_{\infty},\rho_{\infty})+\frac{3}{5}\int_{\mathbb{R}^3}(|x|^{-1}*\rho_{\infty})(\rho_{c}-\rho_{\infty})(x)\,dx\\
=&\frac{3}{10m}\cdot\big(\frac{6\pi^2}{q}\big)^{\frac{2}{3}}\left[ \int_{\mathbb{R}^3}\rho_{c}(x)^\frac{5}{3}\,dx-\int_{\mathbb{R}^3}\rho_{\infty}(x)^\frac{5}{3}\,dx\right]+D(\rho_{c},\rho_{c})-D(\rho_{\infty},\rho_{\infty}).
\end{aligned}
\end{gather}
On the last equality, we use the Euler-Lagrange equation for $\rho_{\infty}$. We split the proof of Theorem \ref{th1.3} (i) in three steps.

\textbf{Step 1:} We claim
\begin{align}\label{eq4.27}
D(\rho_{c},\rho_{c})-D(\rho_{\infty},\rho_{\infty})\ge 0.
\end{align}
By \eqref{eq-Drhoc} and \eqref{eq-Drhoinfty}, we obtain 
\begin{gather}\label{eq4.28}
	\begin{aligned}
		&\kappa D(\rho_{c},\rho_{c})-\kappa D(\rho_{\infty},\rho_{\infty})\\
		&=2\mathcal{E}_{\infty}(\rho_{\infty})-2\mathcal{E}_{c}(\rho_{c})+\int_{\mathbb{R}^3}j_{mc}(\rho_{c}(x))\,dx+m^2c^4\int_{\mathbb{R}^3}\bar{j}_{mc}(\rho_{c}(x))\,dx-mc^2N.
	\end{aligned}
\end{gather}
By the definition of $j_{mc}(\rho_{c})$ and $\bar{j}_{mc}(\rho_{c})$, we find that
\begin{gather}\label{eq4.29}
	\begin{aligned}
		\int_{\mathbb{R}^3}&j_{mc}(\rho_{c}(x))\,dx+m^2c^4\int_{\mathbb{R}^3}\bar{j}_{mc}(\rho_{c}(x))\,dx-mc^2N\\
		=&\frac{q}{8\pi^3}\int_{\mathbb{R}^3}\int_{|p|<\eta_c}\left[(\sqrt{c^{2}|p|^2 +m^2c^{4}}-mc^2)+\frac{m^2c^4}{\sqrt{c^2|p|^2+m^2c^4}}-mc^2\right]\,dpdx\\
		=&\frac{q}{8\pi^3}\int_{\mathbb{R}^3}\int_{|p|<\eta_c}(\sqrt{c^{2}|p|^2 +m^2c^{4}}-mc^2)\,dpdx\\
		&-\frac{q}{8\pi^3}\int_{\mathbb{R}^3}\int_{|p|<\eta_c}\frac{mc^2(\sqrt{c^2|p|^2 +m^2c^4}-mc^2)}{\sqrt{c^2|p|^2 +m^2c^4}}\,dpdx\\
		=&\frac{q}{8\pi^3}\int_{\mathbb{R}^3}\int_{|p|<\eta_c}\left( 1-\frac{mc^2}{\sqrt{c^2|p|^2 +m^2c^4}}\right)(\sqrt{c^2|p|^2 +m^2c^4}-mc^2)\,dpdx\\
		=&\frac{q}{8\pi^3}\int_{\mathbb{R}^3}\int_{|p|<\eta_c}\frac{(\sqrt{c^2|p|^2 +m^2c^4}-mc^2)^2}{\sqrt{c^2|p|^2 +m^2c^4}}\,dpdx\ge0.
	\end{aligned}
\end{gather}
Since $\mathcal{E}_{\infty}(\rho_{\infty})-\mathcal{E}_{c}(\rho_{c})\ge0$, combining  \eqref{eq4.28} and \eqref{eq4.29}, then \eqref{eq4.27} holds. 

\textbf{Step 2: upper bound.}
It follows from \eqref{eq4.26} and \eqref{eq4.27} that
\begin{gather}\label{eq4.30}
\begin{aligned}
&\frac{3}{10m}\cdot\big(\frac{6\pi^2}{q}\big)^{\frac{2}{3}}\left[ \int_{\mathbb{R}^3}\rho_{c}(x)^\frac{5}{3}\,dx-\int_{\mathbb{R}^3}\rho_{\infty}(x)^\frac{5}{3}\,dx\right]\\&\le\mathcal{E}_\infty(\rho_{c})-\mathcal{E}_\infty(\rho_{\infty})\\
&\le\mathcal{E}_\infty(\rho_{c})-\mathcal{E}_c(\rho_{\infty})\\
&\le\mathcal{E}_\infty(\rho_{c})-\mathcal{E}_c(\rho_{c})\lesssim\frac{1}{c^2},\ \ (\text{by Lemma \ref{lem-energy-estimate}}).
\end{aligned}
\end{gather}

 \textbf{Step 3: lower bound.} Note that
\begin{align*}
\int_{\mathbb{R}^3}j_{\infty}(\rho_{c}(x))\,dx-\int_{\mathbb{R}^3}j_{\infty}(\rho_{\infty}(x))\,dx=\frac{3}{5}\cdot\frac{1}{2m}\big(\frac{6\pi^2}{q}\big)^{\frac{2}{3}}\big[\int_{\mathbb{R}^3}\rho_{c}(x)^{\frac{5}{3}}\,dx-\int_{\mathbb{R}^3}\rho_{\infty}(x)^{\frac{5}{3}}\,dx\big].
\end{align*}
Thus, to prove  
\begin{align}
\big[\int_{\mathbb{R}^3}\rho_{c}(x)^{\frac{5}{3}}\,dx-\int_{\mathbb{R}^3}\rho_{\infty}(x)^{\frac{5}{3}}\,dx\big]\gtrsim\frac{1}{c^2}
\end{align}
 is equivalent to prove
\begin{align}\label{eq425}
\int_{\mathbb{R}^3}j_{\infty}(\rho_{c}(x))\,dx-\int_{\mathbb{R}^3}j_{\infty}(\rho_{\infty}(x))\,dx\gtrsim\frac{1}{c^2}.
\end{align}
Note that 
\begin{align*}
j_{\infty}(\rho_{c})-j_{\infty}(\rho_{\infty})= j_{\infty}(\rho_{c})-j_{mc}(\rho_{c})+j_{mc}(\rho_{c})-j_{\infty}(\rho_{\infty}).
\end{align*}
Since
\begin{align}\label{eq427}
\int_{\mathbb{R}^3}[j_{mc}(\rho_{c})-j_{\infty}(\rho_{\infty})]\,dx=\mathcal{E}_{c}(\rho_{c})-\mathcal{E}_{\infty}(\rho_{\infty})+\kappa D(\rho_{c},\rho_{c})-\kappa D(\rho_{\infty},\rho_{\infty}).
\end{align}
Combining \eqref{eq427} and \eqref{eq4.28}, then
\begin{gather}\label{eq429}
\begin{aligned}
&\int_{\mathbb{R}^3}[j_{mc}(\rho_{c})-j_{\infty}(\rho_{\infty})]\,dx\\
&=\mathcal{E}_{\infty}(\rho_{\infty})-\mathcal{E}_{c}(\rho_{c})+\int_{\mathbb{R}^3}j_{mc}(\rho_{c}(x))\,dx+m^2c^4\int_{\mathbb{R}^3}\bar{j}_{mc}(\rho_{c}(x))\,dx-mc^2N.
\end{aligned}
\end{gather}
On the other hand, we find that
\begin{align}\label{eq430}
\int_{\mathbb{R}^3}[j_{\infty}(\rho_{c})-j_{mc}(\rho_{c})]\,dx=\mathcal{E}_{\infty}(\rho_{c})-\mathcal{E}_{c}(\rho_{c}).
\end{align}
Combining \eqref{eq429} and \eqref{eq430}, we have
\begin{gather}\label{eq4.36}
    \begin{aligned}
\int_{\mathbb{R}^3}&j_{\infty}(\rho_{c}(x))\,dx-\int_{\mathbb{R}^3}j_{\infty}(\rho_{\infty}(x))\,dx\\
=&\mathcal{E}_{\infty}(\rho_{c})-\mathcal{E}_{c}(\rho_{c})+\mathcal{E}_{\infty}(\rho_{\infty})-\mathcal{E}_{c}(\rho_{c})\\
&+\int_{\mathbb{R}^3}j_{mc}(\rho_{c}(x))\,dx+m^2c^4\int_{\mathbb{R}^3}\bar{j}_{mc}(\rho_{c}(x))\,dx-mc^2N.
    \end{aligned}
\end{gather}
Since $\int_{\mathbb{R}^3}j_{mc}(\rho_{c}(x))\,dx+m^2c^4\int_{\mathbb{R}^3}\bar{j}_{mc}(\rho_{c}(x))\,dx-mc^2N\ge0$. From \eqref{eq4.36}, we know that
\begin{align*}
	\int_{\mathbb{R}^3}j_{\infty}(\rho_{c}(x))\,dx-\int_{\mathbb{R}^3}j_{\infty}(\rho_{\infty}(x))\,dx
	\ge\mathcal{E}_{\infty}(\rho_{c})-\mathcal{E}_{c}(\rho_{c})+\mathcal{E}_{\infty}(\rho_{\infty})-\mathcal{E}_{c}(\rho_{c}).
\end{align*}
Then \eqref{eq421} and\eqref{eq424} imply 
\begin{align}\label{eq4.37}
	\int_{\mathbb{R}^3}j_{\infty}(\rho_{c}(x))\,dx-\int_{\mathbb{R}^3}j_{\infty}(\rho_{\infty}(x))\,dx\ge2[\mathcal{E}_{\infty}(\rho_{\infty})-\mathcal{E}_{c}(\rho_{\infty})]\gtrsim\frac{1}{c^2}.
\end{align}
Combing \eqref{eq4.36} and \eqref{eq4.37}, we deduce \eqref{eq425} holds. 

Then \eqref{eq4.30} and \eqref{eq4.37} imply Theorem \ref{th1.3} (i) holds.

\subsection{The convergence rate of the compact support}
In this subsection, we consider the convergence rate of the  $R_\infty-R_c$ and complete the proof of the case (ii) of Theorem \ref{th1.3}. 

By Lemma \ref{lem-mul0}, we know that $ \mu_{c}=\frac{\kappa N}{R_{c}}, \mu_{\infty}=\frac{\kappa N}{R_{\infty}}$, thus we need to get the estimate of Lagrange multiplier $\mu_{c}-\mu_{\infty}$. We divide the proof into three steps.

\textbf{Step 1: find the formulation of Lagrange multiplier.} Combining \eqref{eq-Pohozaev-type-identity1} and \eqref{eq-pohozaev-identity2}, we deduce from \eqref{eq-mu_c} and \eqref{eq-mu-infty} that
 \begin{gather}\label{eq437}
 \begin{aligned}
 (\mu_{c}-\mu_{\infty})N=&4\mathcal{E}_{\infty} (\rho_\infty)-4\mathcal{E}_c(\rho_c)+2\int_{\mathbb{R}^3}j_{mc}(\rho_{c}(x))\,dx+2m^2c^4\int_{\mathbb{R}^3}\bar{j}_{mc}(\rho_{c}(x))\,dx\\
 &-2mc^2N+\frac{1}{2m}\int_{\mathbb{R}^3}\eta_{\infty}^2\rho_\infty\,dx-\int_{\mathbb{R}^3}(\sqrt{c^2\eta_{c}^2+m^2c^4}-mc^2)\rho_{c}(x)\,dx.
 \end{aligned}
 \end{gather}
 Note that
 \begin{gather}\label{eq438}
 \begin{aligned}
 &\frac{1}{2m}\int_{\mathbb{R}^3}\eta_{\infty}^2\rho_\infty\,dx-\int_{\mathbb{R}^3}(\sqrt{c^2\eta_{c}^2+m^2c^4}-mc^2)\rho_{c}(x)\,dx\\
 &=\frac{1}{2m}\int_{\mathbb{R}^3}(\eta_{\infty}^2\rho_{\infty}-\eta_{c}^2\rho_{c})\,dx+\int_{\mathbb{R}^3}\left[ \frac{\eta_{c}^2}{2m}-(\sqrt{c^2\eta_{c}^2+m^2c^4}-mc^2)\right]\rho_{c}(x)\,dx.
 \end{aligned}
 \end{gather}
 Direct calculation yields
\begin{gather}
\begin{aligned}
\int_{\mathbb{R}^3}j_{\infty}(\rho_{\infty}(x))\,dx-\int_{\mathbb{R}^3}j_{\infty}(\rho_{c}(x))\,dx&=\frac{3}{5}\cdot\frac{1}{2m}\big(\frac{6\pi}{q}\big)^{\frac{2}{3}}\big[\int_{\mathbb{R}^3}\rho_{\infty}(x)^{\frac{5}{3}}\,dx-\int_{\mathbb{R}^3}\rho_{c}(x)^{\frac{5}{3}}\,dx\big]\\
&=\frac{3}{5}\cdot\frac{1}{2m}\int_{\mathbb{R}^3}(\eta_{\infty}^2\rho_{\infty}-\eta_{c}^2\rho_{c})\,dx.
\end{aligned}
\end{gather}
Then
\begin{gather}\label{eq440}
\begin{aligned}
&\frac{1}{2m}\int_{\mathbb{R}^3}(\eta_{\infty}^2\rho_{\infty}-\eta_{c}^2\rho_{c})\,dx\\
&=\frac{5}{3}\left(\int_{\mathbb{R}^3}j_{\infty}(\rho_{\infty})\,dx-\int_{\mathbb{R}^3}j_{\infty}(\rho_{c})\,dx\right)\\
&=\frac{5}{3}\left( \int_{\mathbb{R}^3}j_{\infty}(\rho_{\infty})\,dx-\int_{\mathbb{R}^3}j_{mc}(\rho_{c})\,dx+\int_{\mathbb{R}^3}j_{mc}(\rho_{c})\,dx-\int_{\mathbb{R}^3}j_{\infty}(\rho_{c})\,dx\right).
\end{aligned}
\end{gather}
The same as \eqref{eq429}, we find that
\begin{gather}\label{eq441}
\begin{aligned}
&\int_{\mathbb{R}^3}j_{\infty}(\rho_{\infty}(x))\,dx-\int_{\mathbb{R}^3}j_{mc}(\rho_{c}(x))\,dx\\
&=\mathcal{E}_{c}(\rho_{c})-\mathcal{E}_{\infty}(\rho_{\infty})-\int_{\mathbb{R}^3}j_{mc}(\rho_{c}(x))\,dx-m^2c^4\int_{\mathbb{R}^3}\bar{j}_{mc}(\rho_{c}(x))\,dx+mc^2N.
\end{aligned}
\end{gather}
By \eqref{eq437}, \eqref{eq438}, \eqref{eq440} and \eqref{eq441}, we obtain
\begin{gather}\label{eq4.43}
\begin{aligned}
(\mu_{c}-\mu_{\infty})N=&\frac{7}{3}(\mathcal{E}_{\infty} (\rho_\infty)-\mathcal{E}_c(\rho_c))+\int_{\mathbb{R}^3}[\frac{\eta_{c}^2}{2m}-(\sqrt{c^2\eta_{c}^2+m^2c^4}-mc^2)]\rho_{c}(x)\,dx.\\
&+\frac{1}{3}\left(\int_{\mathbb{R}^3}j_{mc}(\rho_{c}(x))\,dx+m^2c^4\int_{\mathbb{R}^3}\bar{j}_{mc}(\rho_{c}(x))\,dx-mc^2N\right)\\
&+\frac{5}{3}\int_{\mathbb{R}^3}[j_{mc}(\rho_{c}(x))-j_{\infty}(\rho_{c}(x))]\,dx.
\end{aligned}
\end{gather}

\textbf{Step 2: upper bound.} we claim that
\begin{align}\label{eq4.44}
(\mu_{c}-\mu_\infty)N\lesssim\frac{1}{c^2}.
\end{align}

For the first term in \eqref{eq4.43}, using Lemma \ref{lem-energy-estimate}, we get 
\begin{align}\label{eq4.45}
	\mathcal{E}_{\infty} (\rho_\infty)-\mathcal{E}_c(\rho_c)\lesssim\frac{1}{c^2}.
\end{align}

Next, we estimate the second term. By Lemma \ref{lem-bound-of-operator}, we have
 \begin{align}\label{eta-c-oper1}
 	\frac{\eta_{c}^2}{2m}-(\sqrt{c^2\eta_{c}^2+m^2c^4}-mc^2)\le\frac{\eta_{c}^4}{8m^3c^2}.
 \end{align}
 Thus, we obtain
 \begin{gather}
 \begin{aligned}
 &\int_{\mathbb{R}^3}\left[ \frac{\eta_{c}^2}{2m}-(\sqrt{c^2\eta_{c}^2+m^2c^4}-mc^2)\right]\rho_{c}(x)\,dx\\
 &\le\frac{1}{8m^3c^2}\int_{\mathbb{R}^3}\eta_{c}^4\rho_{c}\,dx\\
 &=\frac{1}{8m^3c^2}\left( \frac{6\pi^2}{q}\right)^{\frac{4}{3}}\int_{\mathbb{R}^3}\rho_{c}(x)^{\frac{7}{3}}\,dx.
 \end{aligned}
 \end{gather}
 Then  Lemma \ref{lem2.3-rhoc-regularity} implies
 \begin{align}\label{eq4.48}
 \int_{\mathbb{R}^3}\left[ \frac{\eta_{c}^2}{2m}-(\sqrt{c^2\eta_{c}^2+m^2c^4}-mc^2)\right]\rho_{c}(x)\,dx\lesssim\frac{1}{c^2}.
 \end{align}

Now, we go to prove 
\begin{align}\label{eq4.49}
	\int_{\mathbb{R}^3}j_{mc}(\rho_{c}(x))\,dx+m^2c^4\int_{\mathbb{R}^3}\bar{j}_{mc}(\rho_{c}(x))\,dx-mc^2N\lesssim\frac{1}{c^2}.
\end{align}
By \eqref{eq4.29}, we have
\begin{gather}
\begin{aligned}
	&\int_{\mathbb{R}^3}j_{mc}(\rho_{c}(x))\,dx+m^2c^4\int_{\mathbb{R}^3}\bar{j}_{mc}(\rho_{c}(x))\,dx-mc^2N\\
	&=\frac{q}{8\pi^3}\int_{\mathbb{R}^3}\int_{|p|<\eta_c}\frac{(\sqrt{c^2|p|^2 +m^2c^4}-mc^2)^2}{\sqrt{c^2|p|^2 +m^2c^4}}\,dpdx.
\end{aligned}
\end{gather}
Using the operator inequality $\sqrt{c^2|p|^2 +m^2c^4}-mc^2\le\frac{|p|^2}{2m}$, we have
\begin{align}
	\frac{(\sqrt{c^2|p|^2 +m^2c^4}-mc^2)^2}{\sqrt{c^2|p|^2 +m^2c^4}}\le\frac{(\frac{|p|^2}{2m})^2}{mc^2}=\frac{|p|^4}{4m^3c^2}.
\end{align}
Then, we have 
\begin{gather}
	\begin{aligned}
		&\int_{\mathbb{R}^3}j_{mc}(\rho_{c}(x))\,dx+m^2c^4\int_{\mathbb{R}^3}\bar{j}_{mc}(\rho_{c}(x))\,dx-mc^2N\\
		&\le\frac{q}{8\pi^3}.\frac{1}{4m^3c^2}\int_{\mathbb{R}^3}\int_{|p|<\eta_c}|p|^4\,dpdx\\
		&=\frac{3}{28m^3c^2}.\left( \frac{6\pi^2}{q}\right)^{\frac{4}{3}}\int_{\mathbb{R}^3}\rho_{c}(x)^\frac{7}{3}\,dx.
	\end{aligned}
\end{gather}
According to Lemma \ref{lem2.3-rhoc-regularity}, then \eqref{eq4.49} holds.

For the last term in \eqref{eq4.43}, using the operator inequality $\sqrt{c^2|p|^2 +m^2c^{4}}\le\frac{|p|^2}{2m}+mc^2$, it's easy to see $$\int_{\mathbb{R}^3}[j_{mc}(\rho_{c}(x))-j_{\infty}(\rho_{c}(x))]\,dx\le 0.$$

Combining \eqref{eq4.45}, \eqref{eq4.48} and \eqref{eq4.49}, we conclude that \eqref{eq4.44} holds.

\textbf{Step 3: lower bound.} We verify 
\begin{align}\label{eq4.53}
(\mu_{c}-\mu_{\infty})N\gtrsim\frac{1}{c^2}.
\end{align}
Since $\frac{1}{2m}\int_{\mathbb{R}^3}\eta_{c}^2\rho_c\,dx=\frac{5}{3}\int_{\mathbb{R}^3}j_{\infty}(\rho_c)\,dx$, then \eqref{eq4.43} becomes
\begin{gather}\label{eq448}
\begin{aligned}
(\mu_{c}-\mu_{\infty})N=&\frac{7}{3}(\mathcal{E}_{\infty} (\rho_\infty)-\mathcal{E}_c(\rho_c))\\
&+\frac{1}{3}\left(\int_{\mathbb{R}^3}j_{mc}(\rho_{c}(x))\,dx+m^2c^4\int_{\mathbb{R}^3}\bar{j}_{mc}(\rho_{c}(x))\,dx-mc^2N\right)\\
&+\frac{5}{3}\int_{\mathbb{R}^3}j_{mc}(\rho_{c}(x))\,dx
-\int_{\mathbb{R}^3}(\sqrt{c^2\eta_{c}^2+m^2c^4}-mc^2)\rho_{c}(x)\,dx.
\end{aligned}
\end{gather}
By calculation, we find that
\begin{gather}\label{eq449}
\begin{aligned}
\frac{5}{3}\int_{\mathbb{R}^3}j_{mc}(\rho_{c}(x))\,dx&=\frac{5}{3}\cdot\frac{q}{8\pi^3}\int_{\mathbb{R}^3}\int_{|p|<\eta_c}(\sqrt{c^{2}|p|^2 +m^2c^{4}}-mc^2)\,dpdx\\
&=\frac{5}{3}\cdot\frac{q}{8\pi^3}\int_{\mathbb{R}^3}\int_{|p|<\eta_c}\frac{c^2|p|^2}{\sqrt{c^{2}|p|^2 +m^2c^{4}}+mc^2}\\
&\ge\frac{5}{3}\cdot\frac{q}{8\pi^3}\int_{\mathbb{R}^3}\int_{|p|<\eta_c}\frac{c^2|p|^2}{\sqrt{c^2\eta_{c}^2 +m^2c^4}+mc^2}\,dpdx\\
&=\int_{\mathbb{R}^3}\frac{c^2\eta_{c}^2}{\sqrt{c^2\eta_{c}^2 +m^2c^4}+mc^2}\rho_{c}(x)\,dx\\
&=\int_{\mathbb{R}^3}(\sqrt{c^2\eta_{c}^2+m^2c^4}-mc^2)\rho_{c}(x)\,dx.
\end{aligned}
\end{gather}
The above inequality implies that
\begin{align}
\frac{5}{3}\int_{\mathbb{R}^3}j_{mc}(\rho_{c}(x))\,dx
-\int_{\mathbb{R}^3}(\sqrt{c^2\eta_{c}^2+m^2c^4}-mc^2)\rho_{c}(x)\,dx\ge 0.
\end{align}
Combining this with $\int_{\mathbb{R}^3}j_{mc}(\rho_{c}(x))\,dx+m^2c^4\int_{\mathbb{R}^3}\bar{j}_{mc}(\rho_{c}(x))\,dx-mc^2N\ge0$, by Lemma \ref{lem-energy-estimate}, we deduce
\begin{align}
(\mu_{c}-\mu_{\infty})N\ge\frac{7}{3}(\mathcal{E}_{\infty} (\rho_\infty)-\mathcal{E}_c(\rho_c))\gtrsim\frac{1}{c^2}.
\end{align}
Thus \eqref{eq4.53} holds.

\textbf{The end proof of Theorem \ref{th1.3} (ii):} Since $\rho_{c}$ and $\rho_{\infty}$ are radial functions, then by Lemma \ref{lem-mul0}, \eqref{eq-Euler-equation1-1} implies that $\rho_c$ has compact support in a ball of radius 
\begin{align}\label{eq-Rc}
	R_{c}=\frac{\kappa N}{\mu_{c}}.
\end{align}
Similarly, we obtain $\rho_{\infty}$ has compact support in a ball of radius 
\begin{align}\label{eq-R-infty}
	R_{\infty}=\frac{\kappa N}{\mu_{\infty}}.
\end{align}
Combining \eqref{eq4.44} and \eqref{eq4.53}, we obtain
\begin{align}\label{eq454}
R_{\infty}-R_{c}=\frac{(\mu_{c}-\mu_{\infty})\kappa N}{\mu_{c}\mu_{\infty}}\sim\frac{1}{c^2}.
\end{align}
This  completes the proof of the Theorem \ref{th1.3} (ii). {{\hfill $\square$}}\\
\subsection{The  decay rates of minimizers}
The purpose of this subsection is to investigate the decay rate of minimizer $\rho_{c}(x)$ and $\rho_{\infty}(x)$, moreover, we complete the proof of Theorem \ref{th1.4}.

 Since $\rho_\infty(x)$ and $\rho_c(x)$ are symmetric decreasing functions, then, by Newton's theorem
\begin{align}
|x|^{-1}*\rho_\infty\leq \frac{N}{|x|}, \ \ |x|^{-1}*\rho_c\leq \frac{N}{|x|}.
\end{align}
Note that $\mu_\infty=\frac{\kappa N}{R_{\infty}}$ and $\mu_c=\frac{\kappa N}{R_c}$, then for $R_{\infty}\ge|x|\geq R_{c}$ we have
	\begin{align}\label{rho-inf-up1}
		\kappa(|x|^{-1}*\rho_\infty)-\mu_{\infty}\leq \frac{\kappa N}{|x|}-\frac{\kappa N}{R_{\infty}}\leq \frac{\kappa N}{R_c}-\frac{\kappa N}{R_{\infty}}=\frac{(R_{\infty}-R_c)\kappa N}{R_cR_{\infty}}\lesssim\frac{1}{c^2},
	\end{align}
and for $R_c-\frac{K_1}{c^2}\leq|x|\leq R_c$
\begin{gather}\label{eq5.16}
		\begin{aligned}
			\kappa(|x|^{-1}*\rho_c)-\mu_c&\leq \frac{\kappa N}{|x|}-\frac{\kappa N}{R_c}\leq\frac{\kappa N}{R_c-\frac{K_1}{c^2}}-\frac{\kappa N}{R_c}\le\frac{2\kappa N}{R_c^2}\frac{K_1}{c^2}.
		\end{aligned}
	\end{gather}
 Since $\eta_{\infty}=(6\pi^2\rho_{\infty}/q)^{\frac{1}{3}}$ and $\rho_\infty(x)=0$ for all $|x|\geq R_c$, then  by \eqref{eq-Euler-equation2-1} and \eqref{rho-inf-up1}, we have
	$$\rho_{\infty}(x)\lesssim  \frac{1}{c^3}.$$
Thus, the case (i) of Theorem \ref{th1.4} holds.	
	
Next, we prove the case (ii) of Theorem \ref{th1.4}. We rewrite \eqref{eq-Euler-equation1-1} as 
	\begin{align*}
		\frac{\eta_{c}^2}{2m}=\big[\frac{\eta_{c}^2}{2m}-(\sqrt{c^2\eta_{c}^2+m^2c^4}-mc^2)]+\kappa\left(|x|^{-1}*\rho_{c}(x)\right)-\mu_{c}, \ \ \text{in $B(R_c)\setminus B(R_c-\frac{K_1}{c^2})$},
	\end{align*}
By \eqref{eta-c-oper1} and lemma \ref{lem2.3-rhoc-regularity} we have
$$\frac{\eta_{c}^2}{2m}-(\sqrt{c^2\eta_{c}^2+m^2c^4}-mc^2)\lesssim\frac{1}{c^2}.$$
Since $\eta_{c}=(6\pi^2\rho_{c}/q)^{\frac{1}{3}}$ and $\rho_c(x)=0$ for all $|x|\geq R_c$, combining \eqref{eq5.16} then 
$$\rho_{c}(x)\lesssim  \frac{1}{c^3}.$$
This completes the proof of the Theorem \ref{th1.4}.{{\hfill $\square$}}

\section{Optimal uniform bounds with respect to $N$}\label{sec5}
Reviewing the previous sections, we observe that the Lagrange multiplier and the compact support of $\rho_c(x)$ and $\rho_\infty(x)$ depending on $N$, therefore, we present some precise estimates with $N$ in this section.

Let $B(R)$ denote a ball centered at the origin with radius $R$.
Note that  
\begin{align}\label{eq5.1}
\sqrt{c^2\eta_{c}^2+m^2c^4}-mc^2=\kappa\left(|x|^{-1}*\rho_{c}(x)\right)-\mu_{c}, \ \ \text{in $B(R_c)$},
\end{align}
and 
\begin{align}\label{eq5.2}
\frac{\eta_{\infty}^2}{2m}=\kappa(|x|^{-1}*\rho_\infty(x))-\mu_{\infty}, \ \ \text{in $B(R_\infty)$},
\end{align}

\begin{lem} \label{lem-5.1} For $m>0$ and $N>0$, we have 
\begin{align}
\mu_{\infty}\sim N^{\frac{4}{3}}, \ \ R_{\infty}\sim N^{-\frac{1}{3}}.
\end{align}
\end{lem}
\begin{pf}
By  Lemma \ref{lem-mul0}, \eqref{eq-mu-infty} and \eqref{eq-pohozaev-identity2}, we have
\begin{align}\label{eq5.4}
\frac{\kappa N}{R_\infty}=\mu_{\infty}=\frac{7}{6N}\kappa D(\rho_{\infty},\rho_{\infty}).
\end{align}
For a fixed $\rho_0(x)\in L^1\cap L^{\frac{5}{3}}$ with $\int_{\mathbb{R}^3}\rho_0(x)\,dx=1$, let $\rho^{\lambda}_N(x)=\lambda^3N\rho_0(\lambda x)$ with $\lambda>0$, then $\int_{\mathbb{R}^3}\rho^{\lambda}_N(x)\,dx=N$, and we have
\begin{gather*}
	\begin{aligned}
		-\kappa D(\rho_\infty, \rho_\infty)&\leq \mathcal{E}_{\infty}(\rho_{\infty})=E_{\infty}(N)\leq \mathcal{E}_{\infty}(\rho^{\lambda}_N)\\
  &= \lambda^2N^{\frac{5}{3}}\int_{\mathbb{R}^3}j_{\infty}(\rho_0)\,dx-\lambda N^2 \kappa D(\rho_0,\rho_0).
	\end{aligned}
\end{gather*}
Take the minimum over $\lambda$ with $\lambda=N^2\kappa D(\rho_0,\rho_0)/(2N^{\frac{5}{3}}\int_{\mathbb{R}^3}j_{\infty}(\rho_0)\,dx)$, then 
\begin{align}
-\kappa D(\rho_\infty, \rho_\infty)\leq -\frac{ N^4\kappa^2 D^2(\rho_0,\rho_0)}{4N^{\frac{5}{3}}\int_{\mathbb{R}^3}j_{\infty}(\rho_0)\,dx} m=-C_0 N^{\frac{7}{3}}.
\end{align}
By \eqref{eq5.4}, then 
\begin{align}\label{eq66}
\mu_{\infty}\geq \frac{7C_0}{6} N^{\frac{4}{3}},\ \ R_{\infty}\leq \frac{6\kappa }{7C_0} N^{-\frac{1}{3}}.
\end{align}

Next, we show the lower bound. From \eqref{eq-mu-infty} and \eqref{eq-pohozaev-identity2}, we obtain
\begin{align}
\mu_{\infty}N\le2\kappa D(\rho_\infty,\rho_\infty)=\frac{6}{5m}\big(\frac{6\pi^2}{q}\big)^{\frac{2}{3}}\int_{\mathbb{R}^3}\rho_\infty^{\frac{5}{3}}(x)\,dx.
\end{align}
Note that \eqref{eq417} implies $\|\rho_{\infty}\|_{L^{\frac{5}{3}}(\mathbb{R}^3)}^\frac{5}{3}\lesssim N^{\frac{7}{3}}$, then we have
\begin{align}
\mu_{\infty}\lesssim N^{\frac{4}{3}},\ \ \frac{\kappa N}{\mu_\infty}=R_{\infty}\gtrsim N^{-\frac{1}{3}}.
\end{align}
Thus, we complete the  proof of this lemma.
\end{pf}

\textbf{The proof of Theorem \ref{th1.5}:} Notice that $R_c\sim N^{-\frac{1}{3}}$ for $c$ large enough  is easy to check by the fact that $R_\infty\sim N^{-\frac{1}{3}}$ and $R_\infty-R_c\sim O(\frac{1}{c^2})$.
Since $\rho_\infty(x),\rho_c(x)$ are symmetric decreasing and continuous, by \eqref{eq5.1} and \eqref{eq5.2} we know that 
 $|x|^{-1}*\rho_{\infty}(x)$ and  $|x|^{-1}*\rho_{c}(x)$ are also symmetric decreasing and continuous. Then
    \begin{align*}
&\|\rho_{\infty}\|_{L^\infty(\mathbb{R}^3)}=\rho_\infty(0),\ \  \|\rho_c\|_{L^\infty(\mathbb{R}^3)}=\rho_c(0),\\
& \|\,|x|^{-1}*\rho_{\infty}\|_{L^{\infty}(\mathbb{R}^3)}=|x|^{-1}*\rho_{\infty}(0),\ \ \ \|\,|x|^{-1}*\rho_{c}\|_{L^{\infty}(\mathbb{R}^3)}=|x|^{-1}*\rho_{c}(0).
     \end{align*}
Since $R_{\infty}\leq N^{-\frac{1}{3}}$, then,  
    \begin{gather}\label{rho-infty0}
 \begin{aligned}
 |x|^{-1}*\rho_{\infty}(0)&=\int_{|y|\leq R_{\infty}}\frac{1}{|y|}\rho_{\infty}(y)\,dy\\
 &\leq \left(\int_{|y|\leq R_\infty}\frac{1}{|y|^\frac{5}{2}}\,dy\right)^{\frac{2}{5}}\left( \int_{|y|\leq R_{\infty}}\rho_{\infty}^{\frac{5}{3}}(y)\,dy\right)^{\frac{3}{5}}\\
 &=(8\pi)^{\frac{2}{5}}R_{\infty}^{\frac{1}{5}}\left( \int_{|y|\leq R_{\infty}}\rho_{\infty}^{\frac{5}{3}}(y)\,dy\right)^{\frac{3}{5}}\\
 &\lesssim N^{-\frac{1}{15}} \rho^{\frac{2}{5}}_\infty(0)   N^{\frac{3}{5}}= \rho^{\frac{2}{5}}_\infty(0)   N^{\frac{8}{15}}.
  \end{aligned}
\end{gather}
Since $\eta_{\infty}=(6\pi^2\rho_{\infty}/q)^{\frac{1}{3}}$ and $\mu_{\infty}>0$, by \eqref{eq5.2} it follows that 
\begin{align}\label{rho-0-up-bound1}
\rho_{\infty}^{\frac{2}{3}}(0)\lesssim \rho_\infty(0)^{\frac{2}{5}} N^{\frac{8}{15}}
\end{align}
This implies that 
\begin{align}
\|\rho_{\infty}\|_{L^\infty}=\rho_{\infty}(0)\lesssim 
N^{2}.
\end{align}
Substitute this into \eqref{rho-infty0}, then
$$|x|^{-1}*\rho_{\infty}(0)\lesssim N^{2\times\frac{2}{5}}N^{\frac{8}{15}}=N^{\frac{4}{3}}.$$

By \eqref{eq5.2} and \eqref{eq66} we have
\begin{align}\label{eq612}
\kappa(|x|^{-1}*|\rho_{\infty}(x))=\frac{\eta_{\infty}^2}{2m}+\mu_{\infty}\ge\mu_{\infty}\gtrsim N^{\frac{4}{3}}.
\end{align}
Combining \eqref{eq612} with \eqref{eq5.2}, it's easy to get
\begin{align}
\rho_\infty(0)\gtrsim N^2.
\end{align}
Thus, case (i) of Theorem \ref{th1.5} holds.

Next, we go to prove case (ii). Since  $R_\infty-R_c=O(\frac{1}{c^2})$ (see Theorem \ref{th1.3}) and $R_{\infty}\lesssim N^{-\frac{1}{3}}$ (see Lemma \ref{lem-5.1}),  then for $c$ large enough we have
$$R_c\lesssim N^{-\frac{1}{3}}.$$
 The similar calculations as \eqref{rho-infty0} and by \eqref{rho-c-53-N},
\begin{align}\label{rho-c0}
 |x|^{-1}*\rho_c(0)\leq(8\pi)^{\frac{2}{5}}R_{c}^{\frac{1}{5}}\left( \int_{|y|\leq R_{c}}\rho_{c}^{\frac{5}{3}}(y)\,dy\right)^{\frac{3}{5}}\lesssim N^{-\frac{1}{15}}N^{\frac{7}{5}}=N^{\frac{4}{3}}.
  \end{align}
The similar to \eqref{eq612}, we can get
\begin{align}\label{eq615}
|x|^{-1}*\rho_c(0)\ge \mu_c\gtrsim N^{\frac{4}{3}}.
\end{align}

On the other hand, take $\delta=2\sqrt{5}m N^{\frac{4}{3}}$ in Lemma \ref{lem2.1-operator-bounded}, for $c$ large enough, we have 
\begin{gather}\label{rho-c1}
 \begin{aligned}
\sqrt{c^2\eta^2_{c}(0)+m^2c^4}-mc^2+2\sqrt{5}m N^{\frac{4}{3}}&\ge B \eta_{c}(0)\\
&=\min\left\lbrace2N^{\frac{2}{3}},\frac{c}{2}\right\rbrace \eta_{c}(0)=2N^{\frac{2}{3}} \eta_{c}(0). 
 \end{aligned}   
\end{gather}
Since $\mu_c>0$, combining \eqref{eq5.1}, \eqref{rho-c0} and \eqref{rho-c1}, then
\begin{align*}2N^{\frac{2}{3}} \eta_{c}(0)&\leq \kappa\left(|x|^{-1}*\rho_{c}(0)\right)-\mu_{c}+ 2\sqrt{5}m N^{\frac{4}{3}} \\
&\leq \kappa\left(|x|^{-1}*\rho_{c}(0)\right)+ 2\sqrt{5}m N^{\frac{4}{3}}
\\
&\lesssim N^{\frac{4}{3}}.
\end{align*}
Since $\eta_c=(6\pi^2\rho_c/q)^{\frac{1}{3}}$, it follows that 
$$\rho^{\frac{1}{3}}_c(0)\lesssim  N^{\frac{2}{3}}.$$
Thus for $c$ large enough
$$\|\rho_c\|_{L^\infty}=\rho_c(0)\lesssim N^{2}.$$
By \eqref{eq5.1} and \eqref{eq615}, we have
\begin{align}
\frac{\eta_c^2}{2m}\ge\sqrt{\eta_{c}^2c^2+m^2c^4}-mc^2=\kappa\left(|x|^{-1}*\rho_{c}(x)\right)-\mu_{c}\gtrsim N^{\frac{4}{3}},
\end{align}
which implies $$\|\rho_c\|_{L^\infty}\gtrsim N^2.$$
This  completes the proof of Theorem \ref{th1.5}.

%%%%%%%%%%%%%%%%%%%%%%%%%%%%%%%%%%%%
\appendix
\section{The existence of minimizers of $E_\infty(N)$}\label{secA}
The existence of minimizers for  limit energy $E_\infty(N)$ is  as following.
\begin{thm}\label{thA1}
	Suppose $q\ge1$ and $m>0$, then $E_\infty(N)$ has at least one minimizer.
\end{thm}

\begin{lem}\label{lem-strictly-binding-inequality}
Fix $q\ge1$ and $m>0$. Then the following strictly binding inequality
\begin{align}\label{eq-A3}
	E_{\infty}(N)<E_{\infty}(\alpha)+ E_{\infty}(N-\alpha)
\end{align}
holds for $0<\alpha<N$.
\end{lem}
\begin{pf}
Firstly, we claim that 
\begin{align}\label{eq-A4}
	E_{\infty}(N)<0.
\end{align}
For any $\lambda>0$, we choose proper scaling $\rho^{\lambda}(x)=\lambda^3\rho(\lambda x)$, with $\int_{\mathbb{R}^3}\rho(x)\,dx=N$, then
\begin{gather}\label{eq-A5}
	\begin{aligned}
		\mathcal{E}_{\infty}(\rho^{\lambda})=\lambda^2\int_{\mathbb{R}^3}j_{\infty}(\rho)\,dx-\lambda \kappa D(\rho,\rho)
	\end{aligned}
\end{gather}
Choose $\lambda$ small enough, then \eqref{eq-A4} holds.

Note that for any $\delta>0$, there exists $v\in L^{\frac{5}{3}}(\mathbb{R}^3)$ with $\|v\|_{L^1}=\mathcal{N}<N$ such that $	E_{\infty}(\mathcal{N})\le\mathcal{E}_{\infty}(v)\le E_{\infty}(\mathcal{N})+\delta.$ Hence for any $\theta>1$ such that $\theta\mathcal{N}\le N$, we have
\begin{gather}\label{eq-A6}
	\begin{aligned}
		E_{\infty}(\theta\mathcal{N})\le&\mathcal{E}_{\infty}(\theta v)=\theta^{\frac{5}{3}}\int_{\mathbb{R}^3}j_{\infty\frac{1}{\theta}}(v(x))\,dx-\theta^2\kappa D(v,v)\\
		=&\theta^2\left(\theta^{-\frac{1}{3}}\int_{\mathbb{R}^3}j_{\infty\frac{1}{\theta}}(v(x))\,dx-\kappa D(v,v)\right)\\
		=&\theta^2\left(\theta^{-\frac{1}{3}}\int_{\mathbb{R}^3}j_{\infty\frac{1}{\theta}}(v(x))\,dx-\int_{\mathbb{R}^3}j_{\infty}(v(x))\,dx\right)+\theta^2\mathcal{E}_{\infty}(v),
	\end{aligned}
\end{gather}
where $j_{\infty\frac{1}{\theta}}(v(x))=\frac{q}{8\pi^3}\int_{|p|<( \frac{6\pi^2v}{q\theta})^{\frac{1}{3}}}\frac{|p|^2}{2m}\,dp$.
Direct calculation, we find that 
\begin{align*}
	\theta^{-\frac{1}{3}}\int_{\mathbb{R}^3}j_{\infty\frac{1}{\theta}}(v(x))\,dx-\int_{\mathbb{R}^3}j_{\infty}(v(x))\,dx=\frac{3}{10m}\left(\theta^{-2}\int_{\mathbb{R}^3}v(x)^{\frac{5}{3}}\,dx-\int_{\mathbb{R}^3}v(x)^{\frac{5}{3}}\,dx\right).
\end{align*}
Since $\theta>1$, thus we deduce
\begin{align}\label{eq-A7}
	\theta^{-\frac{1}{3}}\int_{\mathbb{R}^3}j_{\infty\frac{1}{\theta}}(v(x))\,dx-\int_{\mathbb{R}^3}j_{\infty}(v(x))\,dx<0
\end{align}
Combining \eqref{eq-A6} and \eqref{eq-A7}, we deduce 
\begin{align}\label{eq-A8}
	E_{\infty}(\theta\mathcal{N})\le\theta^2\mathcal{E}_{\infty}(v)\le\theta^2(E_{\infty}(\mathcal{N})+\delta).
\end{align}

Next, we claim that for any $0<\alpha<N$,
\begin{align}\label{eq-A9}
	E_{\infty}(N)<\frac{N}{\alpha}E_{\infty}(\alpha).
\end{align}
Indeed, if $E_{\infty}(\alpha)\ge0$, that \eqref{eq-A9} holds obviously since $E_{\infty}(N)<0$. If $E_{\infty}(\alpha)<0$,
choose $\theta=N/\alpha$ and $\alpha=\mathcal{N}$ in inequality \eqref{eq-A8}, let $\delta<(\theta^{-1}-1)E_{\infty}(\alpha)$, it follows that \eqref{eq-A9} holds. 

In the same way, replace $\alpha$ with  $N-\alpha$, we have
\begin{gather}\label{eq-A10}
	\begin{aligned} 
		E_{\infty}(N)<\frac{N}{N-\alpha}E_{\infty}(N-\alpha)
	\end{aligned}
\end{gather}
Combining \eqref{eq-A9} with \eqref{eq-A10}, then \eqref{eq-A3} holds.
\end{pf}

By Lemma \ref{lem-32bounded}, the minimizing sequence $\{\rho_c\}$ is bounded in $L^{\frac{5}{3}}(\mathbb{R}^3)$. Consequently, there exists a subsequence $\{\rho_{c_k}\}$ such that $\rho_{c_k}\rightharpoonup\rho_{\infty}$. We now apply the
concentration-compactness lemma.
\begin{lem}\label{lem-CMl} 
Let $\{\rho_c\}$ be a bounded sequence in $L^{\frac{5}{3}}(\mathbb{R}^3)$ satisfying $\rho_c\geq0$ and $\int_{\mathbb{R}^3}\rho_c(x)\,dx=N$. Then, there exists a subsequence $\{\rho_{c_k}\}$ satisfying one of the three following possibilities:
\begin{itemize}
    \item[\rm(i)]Compactness: $\rho_{c_k}$ is tight, i.e., for all $\varepsilon>0$, there exists $R<\infty$ such that 
   \begin{align} 
	\int_{|x|\leq R}\rho_{c_k}(x)\,dx\geq N-\varepsilon.
  \end{align}
 \item[\rm(ii)] Vanishing:
	\begin{align} 
	 \lim_{k\to\infty}\int_{|x|\leq R}\rho_{c_k}(x)\,dx= 0,\ \ \text{for all }\  R<\infty.
	\end{align}
\item[\rm(iii)]Dichotomy: There exist $\alpha\in(0, 1)$ and a sequence $\{R_k\}\subset \mathbb{R}^+$ with $R_k\to\infty$ such that
	\begin{gather}
		\begin{aligned} 
		\lim_{k\to\infty}	\int_{|x|\leq R_k}\rho_{c_k}(x)\,dx=\alpha,\qquad\lim_{k\to\infty}	\int_{R_k\leq|x|\leq 6R_k}\rho_{c_k}(x)\,dx=0.
		\end{aligned}
	\end{gather}
		\end{itemize} 
	\end{lem}

Invoking Lemma \ref{lem-CMl}, we obtain a suitable subsequence $\rho_{c_k}$ with $\rho_{c_k}\rightharpoonup\rho_{\infty}$, satisfies either (i), (ii) or (iii). We refer to \cite{JMPNg} for a similar proof of (ii) and\cite{ARMAgf} for a similar proof of (iii).
We rule out (ii) and (iii)  as follows:

\textbf{Vanishing does not occur:} If vanishig occurs, it follows from the similar arguments in \cite{AIHP1Lions1,AIHP1Lions2,CMPjf}  that 
\begin{align*}
\lim_{k\to\infty}D(\rho_{c_k},\rho_{c_k})=0.
\end{align*}
Then, we have  
\begin{align}
    E_{\infty}(N)=\lim_{k\to\infty}\mathcal{E}_{\infty}(\rho_{c_k})=\lim_{k\to\infty}\int_{\mathbb{R}^3}j_{\infty}(\rho_{c_k})\,dx\ge0,
\end{align}
which contradicts $E_{\infty}(N)<0$. Thus, vanishing does not occur.

\textbf{Dichotomy does not occur:} If (iii) is true for $\rho_{c_k}$, the same arguments as Lemma 5 of \cite{JMPNg}, we have 
\begin{align*}
    E_{\infty}(N)\geq E_{\infty}(\alpha)+E_{\infty}(N-\alpha)
\end{align*}
for $0<\alpha<N$. This contradicts the strict binding inequality. Thus, dichotomy  does not occur. Therefore, we have compactness.

\textbf{The end proof of Theorem \ref{thA1}:} From above arguments, we have shown that there exists a subsequence $\rho_{c_k}$ such that (i) of Lemma \ref{lem-CMl} holds true. Then, we have
\begin{align}
    N\ge\int_{\mathbb{R}^3}\rho_{\infty}(x)\,dx\ge\int_{|x|\le R}\rho_{\infty}(x)\,dx=\lim_{k\to\infty}\int_{|x|\le R}\rho_{c_k}(x)\,dx\ge N-\varepsilon
\end{align}
for every $\varepsilon>0$ and suitable $R=R(\varepsilon)<\infty$. This implies that $\int_{\mathbb{R}^3}\rho_{\infty}(x)\,dx=N$. Next, we prove that $\rho_{\infty}$ is a minimizer of $E_{\infty}(N)$. We first deduce from the norm preservation  and the same arguments in   \cite[p.125]{AIHP1Lions1} that we have 
\begin{align}\label{eq:A16}
    \lim_{k\to\infty}D(\rho_{c_k},\rho_{c_k})=D(\rho_{\infty},\rho_{\infty}).
\end{align}
On the other hand, since $\int_{\mathbb{R}^3}j_{\infty}(\rho(x))\,dx$ is weakly lower continuous, then we have 
\begin{align}\label{eq:A17}
    \lim_{k\to\infty}\inf\int_{\mathbb{R}^3}j_{\infty}(\rho_{c_k}(x))\,dx\ge\int_{\mathbb{R}^3}j_{\infty}(\rho_{\infty}(x))\,dx.
\end{align}
Combining \eqref{eq:A16} and \eqref{eq:A17}, we obtain
\begin{align}
    E_{\infty}(N)= \lim_{k\to\infty}\mathcal{E}_{\infty}(\rho_{c_k})\ge\mathcal{E}_{\infty}(\rho_{\infty})\ge E_{\infty}(N),
\end{align}
which implies $\rho_{\infty}$ is a minimizer of $E_{\infty}(N)$. 

This completes the proof of Theorem \ref{thA1}.

\vspace{0.3cm}
\noindent
{\bf Acknowledgments.}

On behalf of all authors, the corresponding author states that there is no conflict of interest. This work has no associated data.

%%%%%%%%%%%%%%%%%%%%%%%%%%%%%%%%%%%%%%%%%%%%%%

\end{document}